\documentclass[12pt,fleqn]{article}
\usepackage{CJK}
\usepackage{latexsym,bm}
\usepackage{amsfonts}
\usepackage{amsthm}
\usepackage{amsmath}
\usepackage{amssymb}
\usepackage{mathrsfs}
\usepackage{graphicx}
\begin{document}

\newcommand{\ep}{\varepsilon}
\newcommand{\p}{\partial}
\newcommand{\be}{\begin{eqnarray}}
\newcommand{\ee}{\end{eqnarray}}
\newcommand{\bee}{\begin{eqnarray*}}
\newcommand{\eee}{\end{eqnarray*}}
\newcommand {\fe} {C^{\alpha,\alpha/2}(\overline\Omega_T)}
\newcommand {\fa} {C^{1+\alpha,(1+\alpha)/2}(\overline\Omega_T)}
\newcommand {\fap} {C^{1+\alpha,(1+\alpha)/2}(\partial\Omega_T)}
\newcommand {\fb} {C^{2+\alpha,1+\alpha/2}(\overline\Omega_T)}
\newcommand {\fbp} {C^{2+\alpha,1+\alpha/2}(\partial\Omega_T)}
\newcommand {\CL}{\cal L}
\newcommand {\CT}{\cal T}
\newcommand {\CU}{\cal U}
\newcommand {\ra}{\rightarrow}

% THEOREMS ---------------------------------------------------------------
\theoremstyle{plain}
\newtheorem{thm}{Theorem}[section]
\newtheorem{cor}[thm]{Corollary}
\newtheorem{lem}[thm]{Lemma}
\newtheorem{prop}[thm]{Proposition}
\newtheorem{exam}[thm]{Example}
\newtheorem{ques}[thm]{Question}
\theoremstyle{definition}
\newtheorem{defn}{Definition}[section]
\newtheorem{ass}{Assumption}[section]
\newtheorem{rmk}{Remark}[section]

\newcommand{\bR}{\mathbb{R}}
\newcommand{\bC}{\mathbb{C}}
\newcommand{\bH}{\mathbb{H}}
\newcommand{\bP}{\mathbb{P}}
\newcommand{\bF}{\mathbb{F}}
\newcommand{\bL}{\mathbb{L}}
\newcommand{\cD}{\mathcal{D}}
\newcommand{\cL}{\mathcal{L}}
\newcommand{\cM}{\mathcal{M}}
\newcommand{\cH}{\mathcal{H}}
\newcommand{\sA}{\mathscr{A}}
\newcommand{\sD}{\mathscr{D}}
\newcommand{\sF}{\mathscr{F}}
\newcommand{\sP}{\mathscr{P}}
\newcommand{\sG}{\mathscr{G}}
\newcommand{\sH}{\mathscr{H}}
\newcommand{\fH}{\mathfrak{H}}

\title{\Large\bf
   Dynkin Game of Stochastic Differential Equations with Random Coefficients,
   and Associated Backward Stochastic Partial Differential Variational Inequality
 \footnote{ \noindent  The work is supported by Basic Research
 Program of China (973 Program, No.2007CB814904), the Science Foundation for Ministry
 of Education of China (No.20090071110001), WCU (World Class University) Program
 through the Korea Science and Engineering Foundation funded by the Ministry of
 Education, Science and Technology (R31-20007), NNSF of China
 (No.10901060, No.10971073), NNSF of GuangDong province (No.9451063101002091),
 Shanghai Postdoctoral Science Foundation (No.10R21411100) and China
 Postdoctoral Science Foundation Funded Project (No.201003244, 20100480555),
 Guangdong InnovativeUniversities in Outstanding Young Talent Cultivation
 Project (No. LYM10062).}}

\author{Zhou Yang~\thanks{School of Mathematical Science, South China Normal
 University, Guangzhou 510631, China.
 \textit{E-mail}:\texttt{yangzhou@scnu.edu.cn} }${\,\,\,}^\ddag$
 ~~~~and~~~~ Shanjian Tang~\thanks{School of Mathematical Science, Fudan
 University, Shanghai 200433, China.
 \textit{E-mail}:\texttt{sjtang@fudan.edu.cn} (Shanjian Tang).}
 ~\thanks{Graduate Department of Financial Engineering, Ajou University,
 San 5, Woncheon-dong, Yeongtong-gu, Suwon, 443-749, Korea.}}

\maketitle
\begin{abstract}
 A  Dynkin game is considered for stochastic differential
 equations with random coefficients. We first apply Qiu and Tang's maximum principle for backward
 stochastic partial differential equations to generalize Krylov
 estimate for the distribution of a Markov process to that of a
 non-Markov process, and establish a generalized It\^o-Kunita-Wentzell's
 formula allowing the test function to be a random field of It\^o's
 type which takes values in a suitable Sobolev space. We then prove the
 verification theorem that the Nash equilibrium point and the value of
 the Dynkin game are  characterized by the strong solution of the
 associated Hamilton-Jacobi-Bellman-Isaacs  equation, which is currently
 a backward stochastic partial differential variational inequality
 (BSPDVI, for short) with two obstacles. We obtain the existence and
 uniqueness result and a comparison theorem for
 strong solution of the BSPDVI. Moreover,
 we study the monotonicity on the strong solution of the BSPDVI by the
 comparison theorem for BSPDVI and define the free boundaries. Finally,
 we identify the counterparts for an optimal stopping time problem
 as a special Dynkin game.
\end{abstract}



\setcounter{equation}{0} \setcounter{section}{0}
 \section{Introduction.}
 Throughout this paper, let $(\Omega,{\cal F},\{{\cal F}_t\}_{t\geq0},\mathbb{P})$
 be a complete filtered probability space, on which we define two independent
 standard Brownian motions: $d_1$-dimension $W\triangleq \{W_t\}_{t\geq0}$
 and  $d_2$-dimension $B\triangleq \{B_t\}_{t\geq0}$. Denote by
 $\mathbb{F}^W\triangleq \{{\cal F}_t^W\}_{t\geq0}$ and
 $\mathbb{F}^B\triangleq \{{\cal F}_t^B\}_{t\geq0}$
 the natural filtrations generated by $W$ and $B$, respectively. We assume that
 they contain all $\mathbb{P}$-null sets in ${\cal F}$. Define
 $\mathbb{F}\triangleq \mathbb{F}^W\vee \mathbb{F}^B$. Denote by
 ${\cal P}$ and ${\cal P}^B$ the $\sigma$- algebras of predictable sets in
 $\Omega\times[\,0,T\,]$ associated with $\mathbb{F}$ and
 $\mathbb{F}^B$, respectively. Denote by ${\cal B}( D)$
 the Borel $\sigma$-algebra  of the domain $ D$ in $\mathbb{R}^{d^*}$.
 \smallskip

 Suppose that the state process $X=(X_1,\cdot\cdot\cdot, X_{d^*})$
 is governed by the following stochastic differential equation
 (SDE, for short):
\be                                                                                              \label{eq1.1}
 X_{i,\,s}=x_i+\int^s_t \beta_{i}(u,X_u)\,du
 +\int^s_t \gamma_{i\,l}(u,X_u)\,dW^l_u+\int^s_t \theta_{i\,k}(u,X_u)\,dB^k_u,
\ee
 where $i=1,\cdot\cdot\cdot,d^*$ and
 $x=(x_1,\cdot\cdot\cdot, x_{d^*})\in \mathbb{R}^{d^*}$.
 The coefficients $\beta, \gamma,$ and $\theta$ are
 ${\cal P}^B\times{\cal B}(\mathbb{R}^{d^*})$
 -measurable random fields taking values in proper spaces,
 and satisfying Assumptions D1 and  D2 (see Section 2 below).
 Note that here and in the following we use repeated indices for summation.
 For example, the repeated subscript $l$ implies summation over
 $l=1,\cdot\cdot\cdot,d_1$ and the repeated subscript $k$ implies
 summation over $k=1,\cdot\cdot\cdot,d_2$. It is clear that SDE
 (\ref{eq1.1}) has a unique strong solution $X^{s,x}$.
 \smallskip

 Let ${\cal U}\,_{t,T}$ be the class of all $\mathbb{F}$-stopping times  which  take
 values in $[\,t,T\,]$. For any
 $(\tau_1,\tau_2)\in {\cal U}\,_{t,T}\times{\cal U}\,_{t,T}$,
 we consider the payoff:
 \bee
 &\!\!\!\!\!\!&R_t\,(x;\tau_1,\tau_2)
 \\[2mm]
 &\!\!\!=\!\!\!&\int_t^{\tau_1\wedge \tau_2}\!\!\!\!\!\!\!\!f_u(X^{t,x}_u)\,du
 +\underline{V}\,_{\tau_1}(X^{t,x}_{\tau_1})\,\chi_{\{\tau_1<\tau_2\wedge T\}}
 +\overline{V}_{\tau_2}(X^{t,x}_{\tau_2})\,\chi_{\{\tau_2\leq\tau_1,\tau_2<T\}}
 +\varphi(X^{t,x}_T)\,\chi_{\{\tau_1\wedge\tau_2= T\}}.\;\qquad
 \eee
 The running cost $f$ and terminal costs $\underline{V}$ and $\overline{V}$ are
 ${\cal P}^B\times{\cal B}(\mathbb{R}^{d^*})$-measurable random fields taking
 values in $\mathbb{R}$, and $\varphi$ is
 ${\cal F}_T^B\times{\cal B}(\mathbb{R}^{d^*})$-measurable random field taking
 values in $\mathbb{R}$. They satisfy Assumptions V3 and V4 (see Section 4 below).
 \medskip

 We consider the following Nash equilibrium of our non-Markovian zero-sum Dynkin
 game (denoted by $\sD_{tx}$ hereafter):  Find a pair
 $(\tau_1^*,\tau_2^*)\in {\cal U}\,_{t,T}\times{\cal U}\,_{t,T}$
 such that for any $\tau_1,\,\tau_2\in{\cal U}\,_{t,T}$, the following
 inequalities hold
 \bee
 \mathbb{E}\,\Big[\,R_t(x;\tau^*_1,\tau_2)\,\Big|\,{\cal F}_t\,\Big]
 \geq \mathbb{E}\,\Big[\,R_t(x;\tau^*_1,\tau_2^*)\,\Big|\,{\cal F}_t\,\Big]
 \geq \mathbb{E}\,\Big[\,R_t(x;\tau_1,\tau_2^*)\,\Big|\,{\cal F}_t\,\Big].
 \eee
 Such a pair $(\tau^*_1,\tau_2^*)$, if it exists,  is called a Nash
 equilibrium point or saddle point of Problem $\sD_{tx}$, and the random variable
 $V_t(x)\stackrel{\triangle}{=}
 \mathbb{E}\,[\,R_t(x;\tau^*_1,\tau_2^*)\,|\,{\cal F}_t\,]$ is
 called the value  of Problem $\sD_{tx}$. We have
 \bee
 V_t(x)=\mathop{{\rm ess.inf}}_{\tau_2\in{\cal U}\,_{t,T}}
 \mathop{{\rm ess.sup}}_{\tau_1\in{\cal U}\,_{t,T}}
 \mathbb{E}\,\Big[\,R_t(x;\tau_1,\tau_2)\,\Big|\,{\cal F}_t\,\Big]
 =\mathop{{\rm ess.sup}}_{\tau_1\in{\cal U}\,_{t,T}}
 \mathop{{\rm ess.inf}}_{\tau_2\in{\cal U}\,_{t,T}}
 \mathbb{E}\,\Big[\,R_t(x;\tau_1,\tau_2)\,\Big|\,{\cal F}_t\,\Big].
 \eee
 The value  $V_t(x)$ is unique if it
 exists. In general, a Nash equilibrium point may not be unique,
 and we shall always mean it by the smallest one.

 Dynkin games were initially introduced by Dynkin and Yushkevich~\cite{Dynkin}, 
 and have received many studies, see among
 others~\cite{Bensoussan,Chang,Cvitanic,G. Peskir,Stettner}.
 Many interesting problems arising from the theory of probability,
 mathematics statistics are reformulated as Dynkin games
 (see \cite{G. Peskir2}). Recently, many new financial problems are formulated
 as Dynkin games, and are turned into partial differential variational
 inequalities (PDVI, for short) or free boundary problems, which are then
 studied via the PDE theory (see~\cite{Friedman1,Yan}, for example).

 The existence of saddle points of Dynkin games has been discussed, either via
 a  pure probabilistic approach, such as Snell's envelope and martingale method,
 or by means of a PDE method. If the coefficients $\beta,\,\gamma,$ and
 $\theta$ of the state equation (\ref{eq1.1}) are all deterministic functions,
 then the state $X$ is Markovian, and the value $V_t(x)$ is a
 deterministic function of $(t,x)$ if the cost functions $f,\,\overline{V},\,
 \underline{V}\,,\,\varphi$ are deterministic function, too. Moreover, it can
 be proved that $V_t(x)$ coincides with the strong solution of the associated
 PDVI by the dynamic programming principle under proper assumptions. See
 Friedman~\cite{Friedman1} for more details. Nowadays very general results on
 a Dynkin game have been established for a right-continuous strong Markov process
 (see Ekstr\"om and Peskir~\cite{EkstroemPeskir} and Peskir~\cite{G. Peskir}).
 In contrast, there are fewer studies on a Dynkin game for a non-Markov process
 (see Lepeltiet and Maingueneau~\cite{LM} and Cvitanic and Karatzas~\cite{Cvitanic}).

 In this paper, we are concerned with the Dynkin game for a non-Markovian process.
 We suppose that the state $X$ is driven
 by two independent standard Brownian motions $W$ and $B$ and the drift
 coefficient $\beta$ and the diffusion coefficients $\gamma,\theta$
 only depend on the path of $B$ in a predictable way, that is,  they are
 all  ${\cal P}^B$-predictable and independent of $W$. The structural
 assumption is used  to guarantee the super-parabolic condition
 (see Assumption V2 in Section 2  below) of the associated backward stochastic
 partial differential variational inequality (BSPDVI, in short). In this
 context,  the value $V_t(x)$ further depends on, in addition to $(t,x)$,
 the path of Brownian motion $B$ up to time $t$. Hence, it is a
 random field. We show that it is characterized by the
 unique strong solution of  the associated Hamilton-Jacobi-Bellman-Isaacs (HJBI)
 equation, which is the following type of BSPDVI:
\be                                                                                                   \label{BSPDI2}
 \left\{
 \begin{array}{l}
 dV_t=-({\cal L} V_t+{\cal M}^k Z^k_t+f_t)\,dt+Z^k_t\,dB^k_t,\qquad
 \mbox{if}\;\;\underline{V}\,_t<V_t<\overline{V}\,_t\,;
 \vspace{2mm} \\
 dV_t\leq-({\cal L} V_t+{\cal M}^k Z^k_t+f_t)\,dt+Z^k_t\,dB^k_t,\qquad
 \mbox{if}\;\;V_t=\underline{V}\,_t\,;
 \vspace{2mm} \\
 dV_t\geq-({\cal L} V_t+{\cal M}^k Z^k_t+f_t)\,dt+Z^k_t\,dB^k_t,\qquad
 \mbox{if}\;\;V_t=\overline{V}\,_t\,;
 \vspace{2mm} \\
 V_T(x)=\varphi(x)\,,
 \end{array}
 \right.
\ee
 where the repeated superscript $k$ is summed from $1$ to $d_2$, and
\be                                                                                                  \label{eq1.3}
 {\cal L}\triangleq a^{ij}\, D_{ij}+b^i\,  D_i+c,\;\;
 {\cal M}^k \triangleq\sigma^{ik}\,  D_i+\mu^k,\;\;
 i,j=1,2,\cdot\cdot\cdot,d^*,\;
 k=1,2,\cdot\cdot\cdot,d_2.\quad
\ee
 The coefficients $a,\,b,\,c,\,\sigma,\,\mu$, the upper obstacle
 $\overline{V}$, and the lower obstacle $\underline{V}\,$ are
 ${\cal P}^B\times{\cal B}(\mathbb{R}^{d^*})$-measurable random fields
 taking values in proper spaces. The terminal value $\varphi$ is
 ${\mathbb{F}}^B_T\times(\mathbb{R}^{d^*})$-measurable random field.

 When the above coefficients are all deterministic, BSPDVI~(\ref{BSPDI2})
 (with the second unknown process $Z$ vanishing) is reduced to a deterministic
 PDVI. There is a huge literature concerning deterministic PDVIs, and see Lions
 and Stampacchia~\cite{LS} and Brezis~\cite{Brezis} among the pioneers and
 Bensoussan and Lions~\cite{BL} and Friedman~\cite{Friedman2} among the monographs.
 On the contrary, there are very few studies on  BSPDVI~(\ref{BSPDI2}) with random
 coefficients. We note that BSPDVIs with one obstacle has been discussed in Chang,
 Pang and Yong~\cite{Chang} in connection with  an optimal stopping problem for an
 SDE with random coefficients, as the associated Hamilton-Jacobi-Bellman
 (HJB, in short) equation, and in  {\O}ksendal, Sulem and Zhang \cite{Zhang}
 in connection with a singular control of SPDEs problem, as a system of backward
 stochastic partial differential equations (BSPDEs, for short)
 with only one reflection, which is a more precise  formulation of BSPDVI with
 only one obstacle. However, they only concern  weak solution of BSPDVIs. 
 In this paper,  we concern strong solution of
 BSPDVI~\eqref{BSPDI2}, which enables us to interpret the derivatives
 $DV,\,D^2V,\,DZ$  almost everywhere in $\Omega\times[\,0,T\,]\times\mathbb{R}^{d^*},$
 and therefore (\ref{BSPDI2}) can be understood point-wisely in $x\in\mathbb{R}^{d^*}$.
 The connection to the associated BSPDVI of the value
 field $V$  is extended to a wider context of the strong solution. The
 existence of such a strong solution requires the super-parabolic condition.

 BSPDVI is in fact a {\it singular} or {\it constrained } BSPDE. It can also be
 regarded as a reflected backward stochastic differential equation (RBSDE, in short)
 in an infinite-dimensional space. Its analysis depends heavily on the state of
 arts of BSPDEs, which is referred to e.g.
 \cite{Bensoussan2,Du,Tang,Karoui,Ma1,Ma2,Peng,Tang1,Tang2,Zhou}. In particular,
 we make use of an estimate by  Du and Tang~\cite{Tang} on the square-integrable
 strong solution theory of BSPDEs in a $C^2$  domain.
 Solution of BSPDVI (\ref{BSPDI2}) is obtained by a  conventional penalty method.
 We consider the penalized approximating BSPDEs~(\ref{eq5.11}), and  show that
 BSPDVI~\eqref{BSPDI2} is  their limit in the strong sense.  The key is to prove the
 convergence  of the nonlinear penalty term.  In~\cite{Bensoussan,Friedman1}, a
 deterministic PDVI is concerned and the convergence is obtained by the compact
 imbedding theorem for Sobolev spaces, which fails to hold in our stochastic
 Sobolev space. In~\cite{Chang},  a BSPDVI with one obstacle is concerned and
 the convergence is obtained via the monotonicity of the approximating BSPDEs' weak
 solution $V_n$ in $n$. Our difficulty has two folds. One comes from the feature
 of two obstacles in our BSPDVI~(\ref{BSPDI2}), which destroys the monotonicity of
 the one-parameterized approximating BSPDEs. The other comes from the strong solution
 of BSPDVI~(\ref{BSPDI2}), which requires an extra  higher order (second-order)
 estimate than that requested for the weak solution.
 We show the convergence by observing that $V_n$ is a Cauchy
 sequence in a proper space, which is an extension of the method
 for RBSDE in~\cite{Cvitanic}.

 To connect Problem $\sD_{tx}$ with BSPDVI~\eqref{BSPDI2}, we have to
 use It\^o-Kunita-Wentzell's formula. The existing one in the literature
 (see Lemma~\ref{lem2.1})  requires that the random test function should
 be twice continuously differentiable. The strong solution of BSPDVI
 (\ref{BSPDI2}) only guarantees that $D^2V$ is  integrable,
 and  not necessarily continuous in general in $x$. Therefore,
 Lemma~\ref{lem2.1} fails to be directly applied to  our computation
 of $V_t(X_t)$, and has to be extended  to more general random test
 functions. We first use Qiu and Tang~\cite{Qiu}'s maximum principle
 for quasilinear BSPDEs to generalize Krylov estimate for the distribution
 of a Markov process to that of a non-Markov
 process. Then using a smoothing method and the generalized Krylov estimate,
 we prove a generalized It\^o-Kunita-Wentzell's formula.

 The rest of the paper is organized as follows. We introduce some notations and results about
 It\^o-Kunita-Wentzell's formula and BSPDEs in Section 2. In Section
 3, we state our hypotheses and generalize It\^o-Kunita-Wentzell's formula. In
 Section 4, we prove the verification theorem that the Nash equilibrium point and
 the value of the Dynkin game are characterized by the strong solution of the
 associated Hamilton-Jacobi-Bellman-Isaacs equation, which is currently a BSPDVI with
 two obstacles. In section 5, we establish the existence and uniqueness result and
 a comparison theorem for strong solution of the BSPDVI, via the strong solution
 theory of BSPDE. In Section 6, we use the comparison theorem for BSPDVI to
 derive properties of the strong solution of BSPDVI (\ref{BSPDI2}),
 and define its stochastic free boundaries under proper assumptions. In the last
 section, we show that the optimal stopping time problem  is a special case of a Dynkin game, and therefore
 similar results hold true here.\medskip

\setcounter{equation}{0} \setcounter{section}{1}
\section{ Preliminaries.}
 In this section, we introduce notations and collect results about
 It\^o-Kunita-Wentzell's formula and BSPDEs.

 Denote by $\mathbb{N}$ and $\mathbb{N}_+$ the set of all nonnegative and positive
 integers, respectively. Denote by $E$ a Euclidean space like $\mathbb{R}$ or
 $\mathbb{R}^{d^*},\,\mathbb{R}^{d^*\times d_1},\,\mathbb{R}^{d^*\times d_2}$, and $
 \mathbb{R}^{d^*\times d^*}$. Moreover, for any
 $x\in \mathbb{R}^{d^*},\,\gamma\in \mathbb{R}^{d^*\times d_1},\,
 \theta\in \mathbb{R}^{d^*\times d_2}$ and $a\in \mathbb{R}^{d^*\times d^*}$, define
 $$
 |\,x\,|\triangleq \left(\sum_{i=1}^{d^*}\,x^2_i\,\right)^{1\over 2}, \quad\quad
 |\,\gamma\,|\triangleq
 \left(\sum_{i=1}^{d^*}\,\sum_{l=1}^{d_1}\,\gamma^2_{il}\,\right)^{1\over 2},
 $$
 $$
 |\,\theta\,|\triangleq
 \left(\sum_{i=1}^{d^*}\,\sum_{k=1}^{d_2}\,\theta^2_{ik}\,\right)^{1\over 2},\quad
 \hbox{ \rm and }|\,a\,|
 \triangleq\left(\sum_{i,\,j=1}^{d^*}\,a^2_{ij}\,\right)^{1\over 2}.
 $$
 Define
 $$
 D_i\triangleq  \p_{x_i};\;\; D_{ij}\triangleq \p_{x_ix_j},\;
 i,j=1,2,\cdots,d^*;\;\;
 D^\alpha \triangleq \p_{x_1}^{\alpha_1}\cdot\cdot\cdot \p_{x_{d^*}}^{\alpha_{d^*}};
 \;\;\mbox{ \rm and }\;\;
 |\,\alpha|\triangleq\sum_{i=1}^{d^*}\alpha_i
 $$
 for any multi-index $\alpha=(\alpha_1,\cdot\cdot\cdot,\alpha_{d^*})$ with
 $\alpha_i\in \mathbb{N}$. Denote by $D\eta$ and $D^2\eta$
 respectively the gradient and the Hessian matrix for a function
 $\eta: E\to \mathbb{R}$.\smallskip

 For an integer $k\in \mathbb{N}$,  $p\in [1,+\infty), q\in [1, +\infty)$, a smooth
 domain $D$ in $\mathbb{R}^{d^*}$, and a positive number $T$, we introduce the
 following spaces:\smallskip

 \noindent$\bullet\,C\,^k\,( D):$ the set of all functions $\eta:D\rightarrow E$ such
 that $\eta$ and $D^\alpha \eta$ are continuous for all

 $\qquad\;\;\;\; 1\leq |\alpha|\leq k$;
 \smallskip

 \noindent$\bullet\,C\,^k_0\,(D):$ the set of  all functions in $C\,^k(D)$ with compact
 support in $D$;
 \smallskip

 \noindent$\bullet\,H^{k,\,p}\,( D):$ the completion of $C\,^k( D)$ under the norm
 $$
 |\,\eta\,|\,_{k,\,p}\triangleq \left(\,\int_ D\,|\,\eta\,|\,^p\,dx
 +\sum_{|\alpha|=1}^k\,\int_ D\,|\,D^\alpha \eta\,|\,^p\,dx\,\right)^{1\over p};
 $$

 \noindent$\bullet\,H_0^{k,\,p}\,( D):$ the completion of $C\,^k_0( D)$ under the norm
 $|\,\eta\,|\,_{k,\,p}\,$;
 \smallskip

 \noindent$\bullet\,\mathbb{L}^{k,p}\,( D):$ the set of all $H^{k,\,p}( D)$-valued and
 $\mathbb{F}^B_T$-measurable random variables such that

 $\qquad\;\;\;\;\;\;\mathbb{E}\,(|\,\varphi|\,_{k,p}^p\,)<\infty$;
 \smallskip

 \noindent$\bullet\,\mathbb{L}_0^{k,p}\,( D):$ the set of all $H_0^{k,\,p}( D)$-valued
 and $\mathbb{F}^B_T$-measurable random variables such that

 $\qquad\;\;\;\;\;\;\mathbb{E}\,(|\,\varphi|\,_{k,p}^p\,)<\infty$;
 \smallskip

 \noindent$\bullet\,{\cal L}\,^p:$ the set of all ${\cal P}$-predictable stochastic
 processes taking values in $E$ with the norm
 $$
 \|X\|_p\triangleq \left[\,\mathbb{E}\left(\int_0^T\,|\,X_t\,|^p\,dt\right)\,
 \right]^{1\over p};
 $$

 \noindent$\bullet\,{\cal S}^p:$ the set of all path continuous processes in
 ${\cal L}\,^p$ with the norm
 $$
 |\|X|\|_p\triangleq \left[\,\mathbb{E}\left(\sup\limits_{t\in[\,0,\,T\,]}
 |\,X_t\,|^p\right) \,\right]^{1\over p};
 $$

\noindent$\bullet\,L_{{\mathbb{F}^B}}^{p,\,q}\,(\mathbb{B}):$ the set of all
 ${\cal P}^B$-predictable stochastic processes with values in Banach

 $\qquad\;\;\;\;$ space $\mathbb{B}$ with the norm
 $$
 \|V\|_{L_{{\mathbb{F}^B}}^{p,\,q}\,(\mathbb{B}\,)}
 \triangleq \left[\,\mathbb{E}\left(\int_0^T\,\|\,V_t\,\|_{\mathbb{B}}^p\,dt\right)
 ^{q\over p}\,\right] ^{1\over q};
 $$

 \noindent$\bullet\,\mathbb{H}^{k,\,p}( D)\triangleq L_{\;{{\mathbb{F}^B}}}^{p,\,p}\,
 (H^{k,\,p}( D))$ with the norm $\|V\|\,_{k,\,p}
 \triangleq \|V\|_{L_{\;{{\mathbb{F}^B}}}^{p,\,p}\,(H^{k,\,p}\,( D\,))}\,$;
 \smallskip

 \noindent$\bullet\,\mathbb{H}\,_0^{k,\,p}( D)\triangleq
 L_{\;{{\mathbb{F}^B}}}^{p,\,p}\,(H_0^{k,\,p}( D))$ with the norm $\|V\|\,_{k,\,p}\,$;
 \smallskip

 \noindent$\bullet\,\mathbb{S}^{k,\,p}( D):$  the set of all path continuous stochastic processes
 in $\mathbb{H}\,^{k,\,p}( D)$ equipped with the norm
 $$
 |\|V|\|\,_{k,\,p}\triangleq
 \left[\,\mathbb{E}\left(\sup\limits_{t\in[\,0,\,T\,]}|\,V_t\,|\,_{k,\,p}^p
 \right)\,\right]^{1\over p};
 $$

 \noindent$\bullet\,\mathbb{M}^p( D):$ the subspace of  $\mathbb{H}\,^{0,\,p}( D)$
 equipped with the norm
 $$
 \| V\|\,_{\mathbb{M}^p}\;\triangleq \;\;
 \mathop{{\rm ess.sup}}_{(w,t)\in\,\Omega\times[\,0,T\,]}
 \left[\, \mathbb{E}\left(\,\int_t^T|\,V_u\,|\,_{0,\,p}^p\;du\,
 \Bigg|\,{\cal F}^B_t\,\right)\,\right]^{1\over p}.
 $$

\begin{rmk}                                                                                         \label{rem2.1}
 The space $\mathbb{M}^p$ is a Banach space. See~\cite{Briand,Qiu}.
\end{rmk}

 The space notations $C\,^k_0\,(\mathbb{R}^{d^*}),\,H^{k,\,p}\,( \mathbb{R}^{d^*})$
 and $\mathbb{L}^{k,p}\,( \mathbb{R}^{d^*}),\,\mathbb{H}^{k,\,p}( \mathbb{R}^{d^*}),\,
 \mathbb{S}^{k,\,p}( \mathbb{R}^{d^*}),\,\mathbb{M}^p( \mathbb{R}^{d^*})$ will be
 abbreviated as $C\,^k_0,\,H^{k,\,p},\,\mathbb{L}^{k,p},\,\mathbb{H}^{k,\,p},\,
 \mathbb{S}^{k,\,p},\,\mathbb{M}^p$ if there is no any confusion.

\begin{rmk}                                                                                         \label{rem2.2}
 We have
 $$
 H^{k,\,p}(\mathbb{R}^{d^*})=H_0^{k,\,p}(\mathbb{R}^{d^*}),\qquad
 \mathbb{L}^{k,\,p}(\mathbb{R}^{d^*})=\mathbb{L}_0^{k,\,p}(\mathbb{R}^{d^*}),\qquad
 \mathbb{H}^{k,\,p}(\mathbb{R}^{d^*})=\mathbb{H}_0^{k,\,p}(\mathbb{R}^{d^*}).
 $$
\end{rmk}

 The following special case of It\^o-Kunita-Wentzell's formula
 (see \cite{Kunita}, \cite{Peng})  is the key to connect Problem $\sD_{tx}$
 and BSPDVI~\eqref{BSPDI2}.

\begin{lem}                                                                                       \label{lem2.1}
 Suppose that the random field
 $V:\Omega\times [\,0,T\,]\times \mathbb{R}^{d^*}\rightarrow  \mathbb{R}$
 satisfies the following:

 (i)$\, V(w,\cdot)$ is continuous with respect to $(t,x)$ a.s. in $\Omega$.

 (ii)$\, V(w,t,\cdot)$ is twice continuously differentiable with respect to $x$ for
 any $t\in[\,0,T\,]$ a.s. in $\Omega$.

 (iii) For each $x\in \mathbb{R}^{d^*},\;V(\cdot,x)$ is a continuous semi-martingale
 of form:
 $$
 V_t(x)=V_0(x)+\int_0^t U_s(x)\,ds+\int_0^t Z^k_s(x)\,dB_s^k
 \;\;\mbox{for any}\;t\in[\,0,T\,]\;\mbox{a.s. in }\Omega,
 $$
 where $U(\cdot,x)$ and $Z(\cdot,x)$ are $\mathbb{F}$-adapted with values in
 $\mathbb{R},\,\mathbb{R}^{d_2}$ for any $x\in \mathbb{R}^{d^*}$, and $Z(w,t,\cdot)$
 is continuously differentiable with respect to $x$ for any $t\in[\,0,T\,]$ a.s. in
 $\Omega$.

 \noindent Let $X$ be a continuous semi-martingale of form (\ref{eq1.1}). Then we have
 \be\nonumber
 V_t(X_t)&=&V_0(X_0)+\int_0^t\,(L V_s+M^k Z^k_s+U_s)\,(X_s)\,ds
 +\int_0^t\,(Z^k_s+M^kV_s)\,(X_s)\,dB^k_s
 \\[2mm]                                                                                              \label{eq2.1}
 &&+\int_0^t\,(N^l\,V_s)\,(X_s)\,dW^l_s,
\ee
 where the repeated superscript $l$ is summed from $1$ to $d_1$ and the repeated
 superscript $k$ is summed from $1$ to $d_2$, and
\be                                                                                                   \label{eq2.2}
 L \triangleq{1\over2}\,(\gamma_{i\,l}\gamma_{j\,l}
 +\theta_{i\,k}\theta_{j\,k})\,D_{ij}
 +\beta_i\,D_i, \;\; M^k \triangleq\theta_{ik}\,D_i,\;\; N^l
 \triangleq\gamma_{il}\,D_i,\;\;
 i,j=1,2,\cdot\cdot\cdot,d^*.\quad
\ee
\end{lem}

 In this paper, we make the following assumptions on the coefficients
 $\beta,\gamma,$ and $\theta$ of SDE (\ref{eq1.1}).
 \smallskip

 \noindent{\bf Assumption D1.} (Boundedness)
 $\beta,\,\gamma,\,\theta$ are ${\cal P}^B\times{\cal B}(\mathbb{R}^{d^*})$
 -measurable with values in $\mathbb{R}^{d^*},\,\mathbb{R}^{d^*\times d_1},
 \,\mathbb{R}^{d^*\times d_2}$, respectively. Moreover, they are bounded by
 a positive constant $K$, i.e.,
 $$
 |\,\beta(\cdot,x)\,|+|\,\gamma(\cdot,x)\,|+|\,\theta(\cdot,x)\,|\leq K\;\;
 \mbox{for any}\;x\in \mathbb{R}^{d^*},\;\mbox{a.e. in}\;\Omega\times(0,T).
 $$

 \noindent{\bf Assumption D2.} (Lipschitz continuity and non-degeneracy)
 $\beta,\,\gamma,$ and $\theta$ satisfy Lipschitz condition in $x$ with constant
 $K$, i.e.,
 \bee
 |\,\beta(\cdot,x_1)-\beta(\cdot,x_2)\,|+|\,\gamma(\cdot,x_1)-\gamma(\cdot,x_2)\,|
 +|\,\theta(\cdot,x_1)-\theta(\cdot,x_2)\,|\leq K\,|x_1-x_2|
 \eee
 for any $x_1,x_2\in \mathbb{R}^{d^*}$, almost everywhere in $\Omega\times(0,T)$.
  Moreover, there exists a positive constant $\kappa$ such that
 $$
 \gamma_{i\,l}\gamma_{j\,l}\,\xi_i\,\xi_j
 \geq \kappa|\,\xi\,|^2\quad
 \mbox{for any}\;\xi\in \mathbb{R}^{d^*},\;
 \mbox{a.e. in}\;\Omega\times(0,T)\times\mathbb{R}^{d^*} .
 $$

 It is clear that SDE (\ref{eq1.1}) has a unique strong solution
 $X\in {\cal S}^p$ for any $p\geq1$.
 \medskip

 Consider the following semilinear BSPDEs  in the domain $D$:
\be  \label{eq2.3}
 \left\{
 \begin{array}{l}
 dV_t(x)=-({\cal L} V_t(x)+{\cal M}^k Z^k_t(x)+F(t,x,V_t(x)))\,dt
 +Z^k_t(x)\,dB^k_t,
 \vspace{2mm} \\
 \hspace{9cm} (t,x)\in[\,0,T)\times D;
 \vspace{2mm} \\
 V_t(x)=0,\;(t,x)\in [\,0,T\,]\times\p  D;
 \qquad V_T(x)=\varphi(x),\;x\in \overline{ D}\;,
 \end{array}
 \right.
\ee
 where $\p  D$ is the boundary of $ D$, and the operators ${\cal L}$ and
 ${\cal M}$ are defined in (\ref{eq1.3}).

\begin{defn}                                                                                         \label{defn2.1}
 (strong solution) The pair of random fields  $(V,Z)\in
 (\mathbb{H}^{2,\,2}( D)\cap\mathbb{H}_0^{1,\,2}( D))\times\mathbb{H}^{1,\,2}( D)$ is called a strong solution of BSPDE~(\ref{eq2.3}),
 if almost surely for all $t\in [\,0,T\,]$, the following holds
 \bee
 V_t(x)=\varphi(x)
 +\int_t^T\,\Big[\,{\cal L} V_s(x)+{\cal M}^k Z^k_s(x)+F(s,x,V_s(x))\,\Big]\,ds
 -\int_t^T\,Z^k_s(x)\,dB^k_s
 \eee
 for a.e. $x\in  D$.
\end{defn}

 Consider the following assumptions on the coefficients in (\ref{eq1.3}).
 \medskip

 \noindent{\bf Assumption V1.} (Boundedness) The given functions
 $a,\,b,\,c,\,\sigma,\,\mu$ are ${\cal P}^B\times{\cal B}(\mathbb{R}^{d^*})$-measurable,
 and are bounded by a constant $K$, taking values  in the set of real symmetric
 $d^*\times d^*$ matrices, in the spaces
 $\mathbb{R}^{d^*},\,\mathbb{R},\,\mathbb{R}^{d^*\times d_2},\,\mathbb{R}^{d_2}$,
 respectively. $Da$ and $D\sigma$ exist almost everywhere and are bounded by $K$.
 \medskip

 \noindent{\bf Assumption V2.} (Super-parabolicity) There exist two positive
 constants $\kappa$ and $K$ such that:
 $$
 \kappa\,|\,\xi|^2+|\sigma^*\xi|^2
 \leq 2\xi^*(a^{ij})\,\xi\leq K\,|\,\xi|^2\quad
 \mbox{for any}\;\xi\in \mathbb{R}^{d^*},\;
 \mbox{a.e. in}\;\Omega\times(0,T)\times\mathbb{R}^{d^*} .
 $$

 In view of~\cite[Theorem~5.3 and Corollary 3.4]{Tang}, we have
 \begin{lem}                                                                                           \label{lem2.2}
 Let Assumptions V1 and V2 be satisfied. Assume that
 $\varphi\in \mathbb{L}_0^{1,2}( D)$ and the random field
 $F:\Omega\times[\,0,T\,]\times D\times\mathbb{R}\rightarrow\mathbb{R}$
 satisfies the following:\smallskip

 (i)$\;F(\cdot,u)$ is ${\cal P}^B\times{\cal B}( D)$-measurable
 for any $u\in \mathbb{R}$;\smallskip

 (ii)$\;F(\cdot,0)\in \mathbb{H}^{0,\,2}( D)$ and $F$
 is Lipschitz continuity with respect to $u$, i.e.,
 $$
 |F(\cdot,v)-F(\cdot,u)|\leq K\,|v-u|\;\;\mbox{for any}\;v,u\in \mathbb{R},\;
 \mbox{a.e. in}\;\Omega\times(0,T)\times D.
 $$
 Then BSPDE~(\ref{eq2.3}) has a unique strong solution $(V,Z)$ such that
 $V\in\mathbb{S}^{1,\,2}( D)$. Moreover, we have the following estimate:
 \bee
 \|V\|^2_{2,\,2}+|\|V|\|^2_{1,\,2}+ \|Z\|^2_{1,\,2} \leq
 C(\kappa,\,K,\,T)\,\Big(\mathbb{E}\,[\,|\,\varphi|^2_{1,\,2}\,]
 +\|\,F(\cdot,V(\cdot))\,\|_{0,\,2}^2\Big).
 \eee
\end{lem}

 According to~\cite[Theorem~5.1]{Tang2}, we give the existence and regularity
 result for strong solution of the linear BSPDE.

\begin{lem}                                                                                            \label{lem2.3}
 Let Assumptions V1 and V2 be satisfied. Assume that $F(\cdot,u)\equiv F(\cdot)$ is
 independent of $u$ and $F\in \mathbb{H}^{k,\,2}( D)$ with $k>d^*/2+2$. Moreover, we
 suppose that there exists a constant $M$ such that the coefficients and terminal
 value satisfy the following
 \bee
 \sum_{1\leq\,|\alpha\,|\,\leq\,k}|\,D^\alpha a\,|
 + |\,D^\alpha b\,|+|\,D^\alpha c\,|
 +|\,D^\alpha \sigma\,|+|\,D^\alpha \mu\,|\leq M,\quad
 \varphi\in \mathbb{L}_0^{1,2}( D)\cap\mathbb{L}^{k+1,2}( D).
 \eee
 Then  BSPDE~(\ref{eq2.3}) has a unique strong
 solution $(V,Z)$ such that $V\in \mathbb{S}^{1,2}( D)
 \cap\,L\,\!_{{\mathbb{F}^B}}^{2,\,2}\,(C\,^2( D')\,)$ and
 $Z\in L\,\!_{{\mathbb{F}^B}}^{2,\,2}\,(C\,^1( D')\,)$ for any domain
 $ D'\subset\subset D$.
\end{lem}

 Now, we recall a special case of the maximum theorem for BSPDE of Qiu and
 Tang~\cite[Theorem~5.8]{Qiu}).

\begin{lem}                                                                                       \label{lem2.4}
 Let Assumptions
 V1 and V2 be satisfied. Assume that $F(\cdot,u)\equiv F(\cdot)$ is independent
 of $u,\,F\in\mathbb{M}^p( D)\cap\mathbb{M}^2( D)$ with $p>d^*+2$, and the
 terminal value satisfies:
 $$
 \mathop{{\rm ess.sup}}_{(w,x)\in\Omega\times D}\,|\,\varphi(w,x)\,|
 <\infty.
 $$
 let $(V,Z)$ be the strong solution  of BSPDE~(\ref{eq2.3}).
 Then we have the following estimate:
 \bee
 \mathop{{\rm ess.sup}}_{(w,t,x)\in\Omega\times[\,0,T\,]\times  D}|V_t(w,x)|
 \leq  C(p,\,\kappa,K,T)\,\left(\|F\|_{\mathbb{M}^p}+\|F\|_{\mathbb{M}^2}
 +\mathop{{\rm ess.sup}}_{(w,x)\in\Omega\times D}\,
 |\,\varphi(w,x)\,|\right).
 \eee
\end{lem}

 Finally, we recall the following backward version It\^o's
 formula for BSPDEs, which is the special case of Theorem 3.2 in \cite{Krylov}:

\begin{lem}                                                                                       \label{lem2.5}
 Let Assumption V1 be satisfied. Assume that
 $$
 (v,Z)\in (\mathbb{H}^{2,\,2}(D)\cap\mathbb{H}_0^{1,\,2}(D))
 \times\mathbb{H}^{1,\,2}(D), \quad
 f\in\mathbb{H}^{0,\,2}(D),\quad \varphi\in \mathbb{L}\,_0^{1,\,2}(D)$$
 and the following equality holds for every $\eta\in C\,_0^2(D)$,
 $$
 \int_D\eta\,v_t\,dx
 = \int_D\eta\,\varphi\,dx
 +\int_t^T\!\!\!\int_D\eta\,({\cal L} v_{s}+{\cal M}^k Z^k_{s}+f_s\,)\,dx\,ds
 -\int_t^T\!\!\!\int_D\eta\,Z^k_{s}\,dx\,dB^k_s
 $$
 almost everywhere in $\Omega\times[\,0,T\,].$
 Then there exists a new version $V$ of $v$ such that
 $V\in\mathbb{S}^{1,\,2}(D)$ and $(V,Z)$ is a strong solution of the following
 BSPDE:
 \bee
 V_t(x)=\varphi(x)
 +\int_t^T\,\Big[\,{\cal L} V_s(x)+{\cal M}^k Z^k_s(x)+f_s(x)\,\Big]\,dt
 -\int_t^T\,Z^k_s(x)\,dB^k_s
 \eee
 with $(t,x)\in [0,T]\times D.$ Moreover, we have the following backward
 version It\^o's formula for BSPDE:

 \bee
 |V_t|_{0,\,2}^2
 &=&|\,\varphi|_{0,\,2}^2+\int_t^T\!\!\!\int_D\,
 \Bigg[\,-2\,a^{ij}  D_iV_s D_jV_s
 +2\left(b^i- D_ja^{ij}\right)  D_iV_sV_s+2c\,V^2_s
 \\[2mm] \nonumber
 &&-2\,\sigma^{ik}Z^k_s  D_iV_s
 +2\left(\mu^k-  D_i\sigma^{ik}\right)\,Z^k_sV_s
 +2f_s\,V_s-\,|Z_s|^2\,\Bigg]\,dx\,ds
 \\[2mm] \nonumber
 &&-2\int_t^T\!\!\!\int_D\,V_s\,Z^k_s\,dx\,dB^k_s\;\;\;\;
 \mbox{for all}\;\;t\in[\,0,T\,]\;\;\mbox{a.s. in}\;\;\Omega.
 \eee
\end{lem}

\setcounter{equation}{0}\setcounter{section}{2}
\section{A generalized It\^o-Kunita-Wentzell's formula.}

 Since a strong solution $V$  is only known to belong
 to $\mathbb{H}^{2,\,2}$, Lemma~\ref{lem2.1} fails to be applied to $V$, and need be
 generalized.  We have

\begin{thm}                                                                                    \label{thm3.1}
 Suppose that the random function
 $V:\Omega\times [\,0,T\,]\times \mathbb{R}^{d^*}\rightarrow\mathbb{R}$
 satisfies the following: $V(x)$ is a continuous semimartigale of form
 $$
 V_t(x)=V_0(x)+\int_0^t U_s(x)\,ds+\int_0^t Z^k_s(x)\,dB_s^k\quad
 \mbox{ a.e.}\;x\in\mathbb{R}^{d^*}
 $$
 for every $t\in[\,0,T\,]$ and almost surely $\omega \in \Omega$,  such that  $V\in \mathbb{H}^{2,\,2},\,Z\in \mathbb{H}^{1,\,2},$ and $
 U\in \mathbb{H}^{0,\,2}$. Let $X$ be a continuous semi-martingale  $X$ of form (\ref{eq1.1}), and
 Assumptions D1 and D2 be satisfied. Then formula (\ref{eq2.1}) still holds.
\end{thm}

\begin{rmk}                                                                                       \label{rem3.1}
 The random function $V$ is not well-defined  at every point $(w,t,x)$ of $\Omega\times [0,T]\times \mathbb{R}^{d^*}$ (since it is defined only in a subset of full measure), and the value of $V_t(w,X_t(w))$ is thus not well-defined in general. However,
 Lemma~\ref{lem3.3} below indicates that $V(X)\triangleq V_{\cdot}(X_{\cdot})$ is
 well-defined as a path-continuous stochastic process in $\Omega\times[0,T]$ (see Remark~\ref{rem3.3} below).
\end{rmk}

\begin{lem}                                                                                           \label{lem3.2}
 Let the assumptions in Lemma~\ref{lem2.3} be satisfied. Assume that
 $F\in\mathbb{M}^p( D)$ with $p>1$, the domain $D$ is bounded, and the terminal value satisfies the following:
 $$
 \mathop{{\rm ess.sup}}_{w\in\Omega} |\,\varphi\,|\,_{0,\,p}<\infty.
 $$
 Let $(V,Z)$ be the unique strong solution of BSPDE~(\ref{eq2.3}). Then we have
 \bee
 \|V\|_{\mathbb{M}^p}\leq  C(p,\,\kappa,K,T)\,\left(\|F\|_{\mathbb{M}^p}+
 \mathop{{\rm ess.sup}}_{w\in\Omega} |\,\varphi\,|\,_{0,\,p}\,\right).
 \eee
 Moreover, if $p\geq 2$, then the above estimate also holds even if $D$ is unbounded.
\end{lem}

\begin{rmk}                                                                                          \label{rem3.2}
Our boundedness assumption on $D$  is used to
 guarantee  $\int_D \Phi_n(0)\,dx <+\infty$ for the case of $p\in (1,2)$ in the following proof. It can be
 removed by using bounded domains to
 approximate $D$ if it is unbounded.
\end{rmk}

 \noindent{\bf Proof of Lemma~\ref{lem3.2}.} If $1<p<2$, we construct the
 mollification of $|\,s|^p$ as
 \bee
 \Phi_n(s)=
 \left\{
 \begin{array}{ll}
 |\,s\,|^p,&{\displaystyle|\,s\,|>{1\over n}};
 \vspace{3mm} \\
 {\displaystyle{p\,(p-2)\over 8\,n^p}\,\,(n^2\,s^2-1)^2
 +{p\over 2n^p}\,(n^2\,s^2-1)+{1\over n^p},}\quad
 &{\displaystyle|\,s\,|\leq{1\over n}},
 \end{array}
 \right.
 \eee
 which uniformly converges to $|\,s\,|^p$. Moreover,
 $\Phi_n\in C^2(\mathbb{R}),\,\Phi'_n(0)=0$, and for any $s\in \mathbb{R}$,
 \bee
 &&\Phi_n(s)\geq0,\qquad  0\leq \Phi''_n(s)\leq C(p,n),\qquad
 |\,\Phi'_n(s)\,|^2\leq C(p)\,\Phi''_n(s)\,\Phi_n(s),
 \\[2mm]
 &&|\,\Phi'_n(s)\,|\leq C(p)\,|\,\Phi_n(s)\,|^{(p-1)/ p},
 \qquad\qquad\quad\,\;\;
 |\,\Phi'_n(s)\,s\,|\leq C(p)\,\Phi_n(s).
 \eee

 Applying It\^o's formula to $\Phi_n(V_t(x))$ for any $x\in  D$ and
 integrating with respect to $x$, we have
 \bee
 &&\int_D\Phi_n(V_t)\,dx\\
 &=&
 \int_t^T\int_D\,\Bigg\{-\Phi''_n(V_s)\,\Big[\,a^{ij}D_iV_s D_jV_s
 +\sigma^{ik}Z^k_s  D_iV_s+{1\over2}\,|Z_s|^2\,\Big]\qquad
 \\[2mm]
 &&+\Phi'_n(V_s)\Big[\,\left(b^i- D_ja^{ij}\right)  D_iV_s
 +c\,V_s+F_s+\left(\mu^k- D_i\sigma^{ik}\right)\,
 Z^k_s\,\Big]\,\Bigg\}\,dx\,ds
 \\[2mm]
 &&-\int_t^T\int_D\Phi'_n(V_s)\,Z^k_s\,dx\,dB^k_s
 +\int_D\Phi_n(V_T)\,dx
 \\[2mm]
 &\leq& \int_t^T\!\!\int_D
 \Bigg\{\Phi''_n(V_s)\,\Bigg[\left(-a^{ij}+\left({1\over2}
 +{\kappa\over8K}\right)\sigma^{ik}\sigma^{jk}\right)
 D_iV_s D_jV_s+{\kappa\over 4}\,|DV_s|^2
 \\[2mm]
 &&+\left({8K+\kappa\over 4(\kappa+4K)}
 -{1\over 2}\right)|Z_s|^2\,\Bigg]
 +C(p,\,\kappa,K)\,(\Phi_n(V_s)+|\,F_s\,|^p)\Bigg\}\,dx\,ds
 \\[2mm]
 &&-\int_t^T\int_D\,\Phi'_n(V_s)\,Z^k_s\,dx\,dB^k_s
 +\int_D\Phi_n(\varphi)\,dx
 \\[2mm]
 &\leq& C(p,\,\kappa,K)\,\,
 \int_t^T\int_D\,(\Phi_n(V_s)+|\,F_s\,|^p)\,dx\,ds
 -\int_t^T\int_D\Phi'_n(V_s)\,Z^k_s\,dx\,dB^k_s
 \\[2mm]
 &&+\int_D\Phi_n(\varphi)\,dx.
 \eee
 Taking the conditional expectation under ${\cal F}^B_t$ on both sides of the above
 inequality and letting $n\rightarrow\infty$,  we deduce that

 \bee
 \int_D |\,V_t\,|^p\,dx\leq
 C(p,\,\kappa,K)\,\mathbb{E}\,\Bigg[\,\int_t^T\int_D\,(\,|\,V_s\,|^p
 +|\,F_s\,|^p\,)\,dx\,ds+|\,\varphi\,|\,^p_{0,\,p}\,\Bigg|\,{\cal F}^B_t\,\Bigg].
 \eee
 Applying the Grownwall inequality, we conclude that
 \bee
 \|V\|^p_{\mathbb{M}^p}\leq T\;\;\;\mathop{{\rm ess.sup}}_{(w,t)\times \Omega\times[\,0,T\,]}
 \int_D |\,V_t\,|^p\,dx
 \leq C(p,\,\kappa,K,T)\,\left(\|F\|^p_{\mathbb{M}^p}+
 \mathop{{\rm ess.sup}}_{w\in\Omega} |\,\varphi\,|\,^p_{0,\,p}\,\right).
 \eee

 If $p\geq2$, we directly apply It\^o's formula to $|\,V_t|^p$. Repeating the
 above argument, we can achieve the desired result.
 \hfill$\Box$\medskip

 To prove Theorem~\ref{thm3.1}, we establish the following priori estimate.

\begin{lem}                                                                                         \label{lem3.3}
 Let $X$ be a continuous semi-martingale of the form (\ref{eq1.1}), and
 Assumptions D1 and D2 be satisfied. If there exist a set $P\in{\cal P}^B$
 and  a domain $ D$ in $\mathbb{R}^{d^*}$ such that
 $f\chi_{P}\chi_D\in \mathbb{M}^2\cap\mathbb{M}^p$, where $\chi_A$ is the indicator
 function of set $A$ and $p>d^*+2$. Then there exists a constant $C$ independent of $f,\,P,\, D$ such that
 \bee
 \mathbb{E}\,
 \Bigg[\,\,\int_0^{T\wedge \tau_P\wedge \tau_D}\!\!\!\!|\,f_t(X_t)|\,dt\,\Bigg]
 \leq C\,(\,\|f\chi_P\chi_D\|_{\mathbb{M}^p}+\|f\chi_P\chi_D\|_{\mathbb{M}^2}\,),
 \eee
 where
 $$
 \tau_P\triangleq \inf\{t\geq0\,:\,\chi_P(t)=0\},\qquad
 \tau_D\triangleq \inf\{t\geq0\,:\,X_t\,\overline{\in}\,D\}.
 $$
If $D$ is bounded and $g\,\chi_{P}\chi_D\in \mathbb{M}^q$ with $1<q\leq p$, then there
 exists a constant $C$ independent of $g,\,P$ such that
 \bee
 \mathbb{E}\,\Bigg[\,\int_0^{T\wedge \tau_P\wedge \tau_ D}
 \left(\int_t^{T\wedge \tau_P\wedge \tau_ D}|\,g_s(X_s)\,|\,ds\right)^{q\over p}
 \,dt\,\Bigg]
 \leq C \,\|g\chi_P\chi_D\|^{q/p}_{\mathbb{M}^q}.
 \eee
 \end{lem}

 \noindent{\bf Proof.} {\it Step 1. }  We prove the first estimate.
 Without loss of generalization, we suppose the random field $f\geq0$.\smallskip

 We shall use It\^o-Kunita-Wentzell's formula and the maximum theorem for
 BSPDE to prove the desired results.  In order to apply It\^o-Kunita-Wentzell's formula,
 we need smooth $f\chi_P\chi_D$.\smallskip

 Define
 \bee
 \widehat{\zeta}(s)&\triangleq&
 \left\{
 \begin{array}{l}
 M\,\exp\left({1\over \,s^2-1}\right)\;\mbox{if}\;|\,s\,|\leq1;
 \vspace{2mm}\\
 0\qquad\qquad\quad\,\;\;\;\mbox{if}\;|\,s\,|>1
 \end{array}
 \right. \eee
 with
 $$M \triangleq\left(\int_{-1}^1 \exp\left({1\over \,s^2-1}\right)\,ds\right)^{-1},\quad
 \zeta(y)=\widehat{\zeta}\left(\,|\,y\,|\,\right),\quad\hbox{ \rm and }\quad
 \zeta_n(y)=n^{d^*}\,\zeta(\,ny\,), $$
and the mollification of $f\chi_P\chi_D$ as
\be                                                                                       \label{eq3.1}
 f^n(w,t,x)\triangleq ((f\chi_P\chi_D)*\zeta_n)(w,t,x)\stackrel{\triangle}{=}
 \int_{\mathbb{R}^{d^*}} (f\chi_P\chi_D)(w,t,y)\,\zeta_n(x-y)\,dy.
\ee
 Then, $f^n\in \mathbb{H}^{k,\,2}$ with $k>d^*/2+2$,
 and there exists a positive constant $C$ such that
 $$
 \| f^n\|_{\mathbb{M}^p}\leq C\,\| f\chi_P\chi_D\|_{\mathbb{M}^p}
 \;\;\hbox{ \rm and}\;\;
 \| f^n\|_{\mathbb{M}^2}\leq C\,\| f\chi_P\chi_D\|_{\mathbb{M}^2}
 \;\;\mbox{for any}\;\;n\in \mathbb{N}_+.
 $$

 Construct $(V^{n,\,m},\,Z^{n,\,m})$ as the solution of the following BSPDE
\be                                                                                       \label{eq3.2}
 \left\{
 \begin{array}{l}
 dV^{n,\,m}_t(x)=-\,(L^m V^{n,\,m}_t(x)+M^m_k Z^{n,\,m}_{k,\,t}(x)+f^n_t(x))\,dt
 +Z^{n,\,m}_{k,\,s}(x)\,dB^k_t,
 \vspace{2mm} \\
 \hspace{9cm} (t,x)\in[\,0,T)\times \mathbb{R}^{d^*};
 \vspace{2mm} \\
 V^{n,\,m}_T(x)=0,\quad x\in \mathbb{R}^{d^*},
 \end{array}
 \right.
\ee
 where
 \bee
 &&L^m \triangleq
 {1\over2}\,(\gamma^m_{i\,l}\,\gamma^m_{j\,l}+\theta^m_{i\,k}\,\theta^m_{j\,k})\,
  D_{ij}+\beta^m_i\,  D_i,\qquad
 M^m_k\triangleq\theta^m_{ik}\, D_i,
 \\[2mm]
 &&\gamma^m_{i\,l}\triangleq\gamma_{i\,l}*\zeta_m,\qquad
 \theta^m_{i\,k}\triangleq\theta_{i\,k}*\zeta_m,\qquad
 \beta^m_i\triangleq\beta_i*\zeta_m.
 \eee

 Since the coefficients of $L$ are bounded and Lipschitz continuous  with respect
 to $x$, then $\gamma^m,\,\theta^m,$ and $\beta^m$ converge to
 $\gamma,\,\theta,\,\beta$ uniformly in $x$, respectively. Moreover, we can check
 that
 \bee
 &&\sum\,_{1\leq|\alpha\,|\leq k}
 \;\;(\,|\,D^\alpha\gamma^m\,|+
 |\,D^\alpha\theta^m\,|+|\,D^\alpha\beta^m\,|\,)
 \leq M_{k,\,m}\;\;\mbox{for any}\;k,\,m\in \mathbb{N}_+;
 \\[2mm]
 &&|\,D\gamma^m\,|+|\,D\theta^m\,|+|\,D\beta^m\,|
 +|\,\gamma^m\,|+|\,\theta^m\,|+|\,\beta^m\,|\leq C\,K;
 \\[2mm]
 &&\gamma^m_{i\,l}\,\gamma^m_{j\,l}\,\xi_i\,\xi_j
 \geq \kappa|\,\xi\,|^2\;\;\mbox{for any}\;\;\xi\in\mathbb{R}^{d^*}.
 \eee

 In view of Lemma~\ref{lem2.3}, BSPDE~(\ref{eq3.2}) admits a strong solution such
 that
 $$
 V^{n,\,m}\in \mathbb{S}^{1,2}(B_j)
 \cap\,L\,\!_{{\mathbb{F}^B}}^{2,\,2}\,(C\,^2(\overline{B}_j)\,)
 \quad\hbox{ \rm and }\quad
 Z^{n,\,m}\in L\,\!_{{\mathbb{F}^B}}^{2,\,2}\,(C\,^1(\overline{B}_j)\,)
 $$
 for any $B_j\triangleq \{x:|\,x|<j\}$ with $j\in \mathbb{N}_+$.\smallskip

 Moreover, the comparison theorem for linear BSPDE (see \cite{Du}) implies that
 $V^{n,\,m}\geq0$. According to Lemma~\ref{lem2.4}, there exists a constant $C$,
 independent of $f,\,P,\,m,\,n$ such that
 \be\nonumber
 \mathop{{\rm ess.sup}}\limits_{w\in\,\Omega}\,
 \sup\limits_{(t,x)\in [0,T]\times D}
 \,V_t^{n,\,m}(x)\,
 \leq C\,(\| f^n\|_{\mathbb{M}^p}+\| f^n\|_{\mathbb{M}^2})
 \leq C\,(\| f\chi_P\chi_D\|_{\mathbb{M}^p}+\| f\chi_P\chi_D\|_{\mathbb{M}^2}).
 \\[2mm]                                                                                                \label{eq3.3}
 \ee

 Let $X^m$ be the strong  solution of the following SDE:
 $$
 X^m_{i,\,s}=x_i+\int^s_t \beta^m_{i}(u,X^m_u)\,du
 +\int^s_t \gamma^m_{i\,l}(u,X^m_u)\,dW^l_u
 +\int^s_t \theta^m_{i\,k}(u,X^m_u)\,dB^k_u.
 $$
 Then $X^m$  converges to $X$ a.e in $\Omega\times[\,0,T\,]$.
 Applying Lemma~\ref{lem2.1} to $V^{n,\,m}_s(X^m_s)$, we have
 \bee
 V^{n,\,m}_0(x)&=&
 \int_0^T(-L^m V^{n,\,m}_s-M^m_k  Z^{n,\,m}_{k,\,s}
 +L^m V^{n,\,m}_s+M^m_k Z^{n,\,m}_{k,\,s}+f^n_s)\,(X^m_s)\,ds\qquad\qquad
 \\[2mm]
 &&-\int_0^T (Z^{n,\,m}_{k,\,s}+M^m_kV^{n,\,m}_s)(X^m_s)\,dB^k_s
 -\int_0^T (\gamma^m_{il}\,D_iV^{n,\,m}_s)(X^m_s)\,dW^l_s.
 \eee
 Taking expectation on both sides of the last equality, in view of
 the estimate~\eqref{eq3.3}, we have
 $$
 \mathbb{E}\,\Bigg[\,\int_0^Tf^n_s(X^m_s)\,ds\,\Bigg]
 =\mathbb{E}\,\Big[\,V^{n,\,m}_0(x)\,\Big]
 \leq C\,(\| f^n\|_{\mathbb{M}^p}+\| f^n\|_{\mathbb{M}^2}).
 $$
 Letting $m\rightarrow +\infty$, we have
 $$
 \mathbb{E}\,\Bigg[\,\int_0^T f^n_s(X_s)\,ds\,\Bigg]
 \leq\liminf\limits_{m\rightarrow\infty}\,
 \mathbb{E}\,\Bigg[\,\int_0^Tf^n_s(X^m_s)\,ds\,\Bigg]
 \leq C\,(\| f^n\|_{\mathbb{M}^p}+\| f^n\|_{\mathbb{M}^2}).
 $$
 Since $\{f^n\}$ is a Cauchy sequence in $\mathbb{M}^p$ and $\mathbb{M}^2$,
 then the above estimate implies that $\{f^n(X)\}$ is a Cauchy sequence in
 ${\cal L}^1$. So, we can think of the definition of $f(X)$ as
 $$
 f(X)
 =\lim\limits_{n\rightarrow\infty} f^n(X)
 \;\;\mbox{in}\;\;{\cal L}^1.
 $$
 Moreover, there exists a subsequence $\{f^{n_k}(X)\}_{k=1}^\infty$
 such that $f^{n_k}(X)$ converges to $f(X)$ a.e. in $\Omega\times[\,0,T\,]$.
 \smallskip

 Since $f\geq0$ and $f^n$ is the mollification of $f\chi_P\chi_D$, then $f^n(X)$
 is nonnegative a.e. in  $\Omega\times[\,0,T\,]$. Hence, we deduce that
 \bee
 \mathbb{E}\,\Bigg[\,\,\int_0^{T\wedge \tau_P\wedge \tau_D}|\,f_t(X_t)|\,dt\,\Bigg]
 &\leq& \mathbb{E}\,\Bigg[\,\,\int_0^T|\,f_t(X_t)|\,dt\,\Bigg]
 =\lim\limits_{n\rightarrow\infty}\mathbb{E}\,\Bigg[\,\,\int_0^T|\,f^n_t(X_t)|\,dt\,\Bigg]
 \\[3mm]
 &\leq&C\,(\| f^n\|_{\mathbb{M}^p}+\| f^n\|_{\mathbb{M}^2})
 \leq  C\,(\| f\chi_P\chi_D\|_{\mathbb{M}^p}+\| f\chi_P\chi_D\|_{\mathbb{M}^2}).
 \eee

 {\it Step 2. } We prove the second estimate. Suppose that $g\geq0$ and smooth $g\chi_P\chi_D$ as
 the above, then $g^n\geq0$. Define $v^{n,\,m}$ as the solution of the following BSPDE:
 \bee
 \left\{
 \begin{array}{l}
 dv^{n,\,m}_t(x)=-\,(L^m v^{n,\,m}_s(x)+M^m_k z^{n,\,m}_{k,\,s}(x)+g^n_s(x))\,ds
 +z^{n,\,m}_{k,\,s}(x)\,dB^k_s,
 \vspace{2mm} \\
 \hspace{9cm} (t,x)\in[\,0,T)\times \mathbb{R}^{d^*};
 \vspace{2mm}\\
 v^{n,\,m}_t(w,x)=0,\;\forall\;(w,t,x)\in \Omega\times \Big(\,\{T\}\times D\cup
 \,[\,0,T\,]\times \p D\,\Big),\qquad
 \end{array}
 \right.
 \eee
 and define $v^{n,\,m}(w,t,x)\equiv0$ if $x\,\overline{\in}\;\overline{D}$.
 \smallskip

 The comparison theorem for linear BSPDE (see \cite{Tang}) implies that
 $v^{n,\,m}\geq0$. Moreover, Lemma~\ref{lem3.2} implies that
 \be                                                                                              \label{eq3.4}
 \|v^{n,\,m}\|_{\mathbb{M}^q}\leq
 C\,\| g^n\,\|_{\mathbb{M}^q}\leq C\,\| g\chi_P\chi_D\|_{\mathbb{M}^q},
 \ee
 where $C$ depends on ${\rm diam}(D)$, and is independent
 of $g,\,P,\,m,\,n$. Hence, we calculate  that

 \bee
 &&\mathbb{E}\,\Bigg[\,\int_0^{T\wedge \tau_P\wedge \tau_D}
 \left(\int_t^{T\wedge \tau_P\wedge \tau_D}\,g^n_s(X^m_s)\,ds\right)
 ^{q\over p}\,dt\,\Bigg]
 \\[2mm]
 &\leq&
 \mathbb{E}\,\Bigg\{\,\int_0^{T\wedge \tau_D}
 \Bigg[\mathbb{E}\,\Bigg(\int_t^{(T\wedge \tau_D)\vee t}\,g^n_s(X^m_s)\,ds\,
 \Bigg|\,{\cal F}_t\,\Bigg)\,\Bigg]^{q\over p}\,dt\,\Bigg\}
 \\[2mm]
 &=&\mathbb{E}\,\Bigg[\,\int_0^{T\wedge \tau_D}
 \Bigg(v^{n,\,m}_t(X^m_t)\Bigg)^{q\over p}\,dt\,\Bigg]
 \hspace{2cm}\mbox{( by Lemma~\ref{lem2.1} and the above method)}
 \\[2mm]
 &\leq& C\,(\,\|(v^{n,\,m})^{q/p}\|_{\mathbb{M}^p}
 +\|(v^{n,\,m})^{q/p}\|_{\mathbb{M}^2}\,)
 \hspace{3.6cm}\mbox{\rm (by the result in Step 1)}
 \\[2mm]
 &\leq& C\|(v^{n,\,m})^{q/p}\|_{\mathbb{M}^p}
 =C\|v^{n,\,m}\|^{q/p}_{\mathbb{M}^q}
 \leq C\,\| g\chi_P\chi_D\|^{q/p}_{\mathbb{M}^q}.
 \hspace{1.3cm}\mbox{\rm (by the estimate~\eqref{eq3.4})}
 \eee

 First, taking $m\rightarrow +\infty$ and then letting $n\rightarrow +\infty$, we
 deduce that
 \bee
 &&\mathbb{E}\,\Bigg[\,\int_0^{T\wedge \tau_P\wedge \tau_D}
 \left(\int_t^{T\wedge \tau_P\wedge \tau_D}|\,g_s(X_s)\,|\,ds
 \right)^{q\over p}\,dt\,\Bigg]
 \\[2mm]
 &\leq&\liminf\limits_{n\rightarrow\infty}\liminf\limits_{m\rightarrow\infty}
 \mathbb{E}\,\Bigg[\,\int_0^{T\wedge \tau_P\wedge \tau_D}
 \left(\int_t^{T\wedge \tau_P\wedge \tau_D}
 |\,g^n_s(X^m_s)\,|\,ds\right)^{q\over p}\,dt\,\Bigg]
 \leq  C \|g\chi_P\chi_D\|^{q/p}_{\mathbb{M}^q}.\qquad
 \eee
 \hfill$\Box$

\begin{rmk}                                                                                \label{rem3.3}
 For any $V\in \mathbb{H}^{0,\,2}$, define
 $$
 P_m\triangleq \left\{(w,t):\;
 \mathbb{E}\left(\,\int_0^T|\,V_u\,|\,_{0,\,2}^2\;du\Bigg|\,{\cal F}_t\,\right)
 \,<m\,\right\}\quad\hbox{ \rm and }\quad
 v\triangleq {\rm sign}(V)\,|\,V\,|^{2/(d^*+3)}.
 $$
 Then $v\chi_{P_m}\chi_{B_m}\in \mathbb{M}^2\cap\mathbb{M}^{d^*+3}$ for any
 $m\in \mathbb{N}_+$. Applying the first estimate in Lemma~\ref{lem3.3} , we
 deduce that $v(X)\,\chi_{\{t\,\leq\,T\wedge\tau_{P_m}\wedge\tau_{B_m}\}}
 \in{\cal L}^1$ and $v(X)$ is well-defined in the set
 $\{(\omega,t):t\,\leq\,T\wedge\tau_{P_m}(\omega) \wedge\tau_{B_m}(\omega)\}$
 for any $m\in \mathbb{N}_+$.  Since
 $P_m\times B_m\uparrow \Omega\times[\,0,T\,]\times \mathbb{R}^{d^*}$,
 then $T\wedge\tau_{P_m}\wedge\tau_{B_m}\uparrow T$ and $v(X)$ is well-defined
 in $\Omega\times[\,0,T\,]$. Hence, the process $V(X)\triangleq {\rm sign}(v(X))\,|\,v(X)\,|^{(d^*+3)/2}$ is well-defined.
\end{rmk}

 \noindent{\bf Proof of Theorem~\ref{thm3.1}.} Smooth $V,\,U,$ and $Z$ as follows
 $$V^n\triangleq V*\zeta_n,\qquad
 U^n\triangleq U*\zeta_n,\qquad \,Z^n\triangleq Z*\zeta_n.$$
 Then we have that for any $x\in \mathbb{R}^{d^*}$,
 $$
 V^n_t(x)=V^n_0(x)+\int_0^t U^n_s(x)\,ds+\int_0^t Z^{n,\,k}_s(x)\,dB_s^k
 \;\;\mbox{for any}\;t\in[\,0,T\,],\;\mbox{a.s. in }\Omega,
 $$
 and $V^n,\,U^n,$ and $Z^n$ satisfy the assumptions in Lemma~\ref{lem2.1} and
 converge to $V,\,U,$ and $Z$ in the spaces
 $\mathbb{H}^{2,\,2},\,\mathbb{H}^{0,\,2}$,
 and $\,\mathbb{H}^{1,\,2}$, respectively.

 Denote $\tau_m\triangleq \tau_{P_m}\wedge\tau_{B_m}$, with $P_m$ waiting to
 be defined in (\ref{eq3.6}). Applying Lemma~\ref{lem2.1}, we deduce that
\be\nonumber
 V^n_{t\wedge\tau_m}(X_{t\wedge\tau_m})
 &=&V^n_{T\wedge\tau_m}(X_{T\wedge\tau_m})
 -\int_{t\wedge\tau_m}^{T\wedge\tau_m}(L V^n_s
 +M^k Z^n_{k,\,s}+U^n_s)\,(X_s)\,ds
 \\[2mm]                                                                                        \label{eq3.5}
 &&-\int_{t\wedge\tau_m}^{T\wedge\tau_m}(Z^n_{k,\,s}+M^kV^n_s)(X_s)\,dB^k_s
 -\int_{t\wedge\tau_m}^{T\wedge\tau_m}(N^l V^n_s)(X_s)\,dW^l_s.\qquad\qquad
\ee

 First, we prove that there is a subsequence (still denoted by itself)
 such that
 $$
 \lim_{n\to \infty}\int_{t\wedge\tau_m}^{T\wedge\tau_m} f^n_s(X_s)\,ds=
 \int_{t\wedge\tau_m}^{T\wedge\tau_m} f_s(X_s)\,ds$$
 and $$
\lim_{n\to \infty} \int_{t\wedge\tau_m}^{T\wedge\tau_m} g^n_{k,\,s}(X_s)\,dB^k_s
 =
 \int_{t\wedge\tau_m}^{T\wedge\tau_m} g_{k,\,s}(X_s)\,dB^k_s
 $$
 almost everywhere in $\Omega\times(0,T)$, where
 $$
 f^n\triangleq L V^n+M^k Z^n_k+U^n,\;\;f\triangleq L V+M^k Z_k+U,\;\;
 g^n_k\triangleq Z^n_k+M^kV^n,\;\;g_k=Z_k+M^kV.
 $$

 Since $V^n$ converges to $V$ in $\mathbb{H}^{2,\,2},\,Z^n$ converges to $Z$
 in $\mathbb{H}^{1,\,2}$ and $U^n$ converges to $U$ in $\mathbb{H}^{0,\,2}$,
 then $f^n$ converges to $f$ in $\mathbb{H}^{0,\,2}$ and
 $g_n$ converges to $g$ in $\mathbb{H}^{1,\,2}$. The Sobolev imbedding
 theorem implies that $g_n$ converges to $g$ also in
 $L\,_{{\mathbb{F}^B}}^{2,\,2}\,(H^{0,\,2q}\,)$,
 where $q=2$ if $d^*\leq 2$ and $q={d^*\over d^*-2}$ if $d^*>2$.  Hence, we
 have that
 $$
 |\,f^n-f\,|\,^2_{0,\,2}+|\,(g^n-g)^2\,|\,_{0,\,q}\rightarrow0\;\;
 \mbox{in}\;\;\mathcal{L}^1.
 $$
 So, there exists a strictly increasing
 sequence of numbers $\{K_m\}_{m=1}^\infty$ such that
 $$
 \mathbb{Q}\left\{(w,t)\in \Omega\times[\,0,T\,]:\,
 |\,f^n-f\,|\,_{0,\,2}+|\,(g^n-g)^2\,|\,_{0,\,q}\geq {1\over m}\right\}
 <{1\over 2^m}\;\;\mbox{for any}\;n\geq K_m,
 $$
 where $\mathbb{Q}$ is the product measure of $\mathbb{P}$ and Lebesgue measure
 on $[\,0,T\,]$.\smallskip

 Define
\be                                                                                                \label{eq3.6}
 P_m\triangleq \bigcap\limits_{n=m}^\infty\left\{(w,t)\in  \Omega\times[\,0,T\,]:
 \,\Big|\,f^{K_n}-f\,\Big|_{0,\,2}+\Big|\,\left(g^{K_n}-g\right)^2\,\Big|_{0,\,q}
 <{1\over n}\,\right\}.
\ee
 We have that
 \bee
 &&\mathbb{Q}(P_m)\geq 1- {1\over 2^{m-1}}\rightarrow 1\;\;
 \mbox{as}\;\;m\rightarrow\infty,
 \\[2mm]
 &&\Big|\,(f^{K_n}-f)\,\chi_{P_m}\chi_{B_m}\,\Big|_{0,\,2}
 +\Big|\,\left(g^{K_n}-g\right)^2\,\chi_{P_m}\chi_{B_m}\,\Big|_{0,\,q}
 \leq {1\over n}\;\;\mbox{for any}\;\;n\geq m,
 \\[2mm]
 &&\|\,(f^{K_n}-f)\,\chi_{P_m}\chi_{B_m} \|_{\mathbb{M}^2}
 +\|\,\left(g^{K_n}-g\right)^2\,\chi_{P_m}\chi_{B_m}\|_{\mathbb{M}^q}
 \rightarrow 0\;\;\mbox{as}\;\;n\rightarrow\infty\;\;
 \mbox{for any}\;\;m\in\mathbb{N}_+.\qquad
 \eee
 From Lemma~\ref{lem3.3} , we obtain that as $n\rightarrow\infty$,
 \bee
 \mathbb{E}\,\Bigg[\,\int_0^{T\wedge\tau_m}
 \left|\int_t^{T\wedge\tau_m}
 \left(f^{K_n}_s(X_s)-f_s(X_s)\right)\,ds\right|^{2\over p}\,dt\,\Bigg]
 \leq C\|\,(f^{K_n}-f)\,\chi_{P_m}\chi_{B_m}\|^{2/ p}_{\mathbb{M}^2}\rightarrow 0.
 \eee
 So, we deduce that there exists a subsequence (still denoted by itself) such that,
 $$
 \int_{t\wedge\tau_m}^{T\wedge\tau_m} f^n_s(X_s)\,ds\rightarrow
 \int_{t\wedge\tau_m}^{T\wedge\tau_m} f_s(X_s)\,ds\;\;
 \mbox{a.e. in}\;\;\Omega\times(0,T)\;\;
 \mbox{as}\;\;n\rightarrow \infty.
 $$
 Moreover,  as $n\rightarrow\infty$,
 \bee
 &&\mathbb{E}\,\Bigg[\,\int_0^{T\wedge\tau_m}\,
 \left|\,\int_t^{T\wedge\tau_m}\,
 \left(\,g^{K_n}_{k,\,s}(X_s)-g_{k,\,s}(X_s)\,\right)\,
 dB^k_s\,\right|^{2q\over p}\,dt\,\Bigg]
 \\[2mm]
 &=&\mathbb{E}\,\Bigg\{\,\int_0^{T\wedge\tau_m}\mathbb{E}\,\Bigg[\,
 \left|\int_t^{T\wedge\tau_m} \left(\,g^{K_n}_{k,\,s}(X_s)-g_{k,\,s}(X_s)\,\right)
 \,dB^k_s\right|^{2q\over p}\,\Bigg|\,{\cal F}_t\,\Bigg]\,dt\,\Bigg\}
 \\[2mm]
 &\leq&C\,\mathbb{E}\,\Bigg\{\,\int_0^{T\wedge\tau_m}\mathbb{E}\,\Bigg[\,
 \left(\int_t^{T\wedge\tau_m} \Big|\,g^{K_n}_s(X_s)-g_s(X_s)\,\Big|^2\,ds\right)
 ^{q\over p}\,\Bigg|\,{\cal F}_t\,\Bigg]\,dt\,\Bigg\}
 \\[2mm]
 &=&C\,\mathbb{E}\,\Bigg[\,\int_0^{T\wedge\tau_m}
 \left(\int_t^{T\wedge\tau_m} \Big|\,g^{K_n}_s(X_s)-g_s(X_s)\,\Big|^2\,ds\right)
 ^{q\over p}\,dt\,\Bigg]
 \\[2mm]
 &\leq& C\left\|\,\left(g^{K_n}-g\right)^2\,\chi_{P_m}\chi_{B_m}\right\|^{q/p}_{\mathbb{M}^q}
 \rightarrow 0.
 \eee

 For any fixed $m$, we pass to the limit in (\ref{eq3.5}) as $n\rightarrow\infty$
 (at least for a subsequence). Using the above method, we easily achieve that a.e.
 in $\{(w,t):0\leq t\leq \tau_m(w),w\in \Omega\}$,
 \bee
 V_{t\wedge\tau_m}(X_{t\wedge\tau_m})&=&V_{T\wedge\tau_m}(X_{T\wedge\tau_m})
 -\int_{t\wedge\tau_m}^{T\wedge\tau_m}(L V_s+M^k Z_{k,\,s}+U_s)\,(X_s)\,ds
 \\[2mm]
 &-&\int_{t\wedge\tau_m}^{T\wedge\tau_m} (Z_{k,\,s}+M^kV_s)(X_s)\,dB^k_s
 -\int_{t\wedge\tau_m}^{T\wedge\tau_m} (N^l\,V_s)(X_s)\,dW^l_s.\qquad
 \eee
 Hence, we can choose a  version of $V(X)$ such that $V_t(X_t)(w)$ is
 continuous with respect to $t$ on $[\,0,\tau_m(w)\,]$ a.s. $w\in\Omega$, and
 the above equality holds for all $t\in [\,0,\tau_m(w)\,]$  a.s. $w\in\Omega$.

 So, we deduce that a.s. in $\Omega$, for all $0\leq t\leq \tau_m$,
 \bee
 V_{t\wedge\tau_m}(X_{t\wedge\tau_m})&=&V_0(X_0)
 +\int_0^{t\wedge\tau_m}(L V_s+M^k Z_{k,\,s}+U_s)\,(X_s)\,ds
 \\[2mm]
 &+&\int_0^{t\wedge\tau_m} (Z_{k,\,s}+M^kV_s)(X_s)\,dB^k_s
 +\int_0^{t\wedge\tau_m} (N^l\,V_s)(X_s)\,dW^l_s.
 \eee
 Taking $m\rightarrow\infty$, we have the desired result.
 \hfill$\Box$

\setcounter{equation}{0} \setcounter{section}{3}
\section{Verification theorem.}

 In this section, we prove the verification theorem that the Nash equilibrium point
 and the value of the Dynkin game are characterized by the strong solution
 of BSPDVI (\ref{BSPDI2}). \smallskip

 Consider the following assumptions on the free term $f$, the terminal value
 $\varphi$, and the upper and lower obstacles  $\overline{V}$ and $\underline{V}\,$ in BSPDVI~\eqref{BSPDI2}.
 \medskip

 \noindent{\bf Assumption V3.} (Regularity)
 $\;f\in \mathbb{H}^{0,\,2},\, \varphi\in \mathbb{L}^{1,\,2}$,  and
 $\underline{V}$ and $\overline{V}$ are continuous semimartigales
 of the following form
 \be                                                                                \label{specialform}
 d\underline{V}\,_t=-\underline{g}\,_t\,dt+\underline{Z}\,^k_t\,dB_t^k,\qquad
 d\overline{V}_t=-\overline{g}_t\,dt+\overline{Z}\,^k_t\,dB_t^k,
 \ee
 where $\underline{g}\,,\,\overline{g}\,,\,
 \underline{Z}\,,\,\overline{Z}\,\in \mathbb{H}^{0,\,2},\,
 \underline{V}\,,\,\overline{V}\,\,\in \mathbb{H}^{1,\,2}$
 and there exists a nonnegative random field $h\in \mathbb{H}^{0,\,2}$ such that
 $$
 {\cal L}\,\underline{V}\,+{\cal M}^k \underline{Z}\,^k-\,\underline{g}\,
 +h \geq 0\quad \hbox{ \rm and } \quad
 {\cal L}\,\overline{V}\,+{\cal M}^k \overline{Z}\,^k-\,\overline{g}\,
 -h \leq 0$$
 hold in the sense of  distribution,
 that is, for any nonnegative function $\eta\in C\,^2_0(\,\mathbb{R}^{d^*})$,
 we have
 \be                                                                                    \label{eq4.1}
 {\cal T}(\,\underline{V}\,,\,\underline{Z}\,,\,
 \underline{g}\,,\,h\,,\,\eta\,) \geq0 \quad \hbox{ \rm and } \quad
 {\cal T}(\,\overline{V}\,,\,\overline{Z}\,,\,
 \overline{g}\,,\,-h\,,\,\eta\,) \leq0\;\;
 \mbox{a.e. in}\;\;\Omega\times[\,0,T\,], \quad
 \ee
 where
 \bee
 {\cal T}(\,\underline{V}\,,\,\underline{Z}\,,\,\underline{g}\,,\,
 h\,,\,\eta\,)
 &\triangleq&-\int_{\mathbb{R}^{d^*}}\Big(\,a^{ij}\,D_i \,\underline{V}\,
 +\sigma^{jk}\,\underline{Z}\,^k\,\Big)\,D_j\eta\,dx
 +\int_{\mathbb{R}^{d^*}}(\,h-\,\underline{g}\;)\,\eta\,dx
 \\[2mm]
 &&+\int_{\mathbb{R}^{d^*}}\Big[\,\left(\,b^i- D_ja^{ij}\,\right)
  D_i\,\underline{V}\,+c\,\underline{V}\,
 +\left(\,\mu^k-  D_i\sigma^{ik}\,\right)
 \,\underline{Z}\,^k\,\Big]\,\eta\,dx.\qquad
 \eee

 \noindent{\bf Assumption V4.} (Compatibility)
 $\;\underline{V}\leq \overline{V},\,\underline{V}\,_T
 \leq \varphi\leq \overline{V}_T$
 and
 $$
 \Big(\,{\cal L}\,\underline{V}\,_t+{\cal M}^k\,\underline{Z}\,^k_t
 +f_t-\,\underline{g}\,_t\,\Big)\;\chi_{\{\,\underline{V}= \overline{V}\,\}}=0
 \;\;\mbox{a.e. in}\;\;\Omega\times \overline{Q},\;\;
 \mbox{where}\;\;Q\triangleq [\,0,T)\times \mathbb{R}^{d^*}.
 $$

 The following is a stronger version of Assumption V$3$.
 \medskip

 \noindent{\bf Assumption V$3'$.}
 $\;f\in \mathbb{H}^{0,\,2},\,\varphi\in \mathbb{L}^{1,\,2}$, and
 $\underline{V}$ and $\overline{V}$ have the following representation:
 $$
 d\underline{V}\,_t=-\underline{g}\,_t\,dt+\underline{Z}\,^k_t\,dB_t^k,\qquad
 d\overline{V}_t=-\overline{g}_t\,dt+\overline{Z}\,^k_t\,dB_t^k,
 $$
 with $\underline{V}\,,\,\overline{V}\,\in \mathbb{H}^{2,\,2},\,
 \underline{Z}\,,\,\overline{Z}\,\in \mathbb{H}^{1,\,2},$ and $
 \underline{g}\,,\,\overline{g}\,\in \mathbb{H}^{0,\,2}$.
 \medskip

 The following clarifies the relationship between Assumptions V$3$ and  V$3'$.

 \begin{prop}                                                                                  \label{prop4.1}
(i)  Assumption V$3'$ implies Assumption V$3$.
(ii) If $\underline{V}$ and $\overline{V}$ satisfy
 Assumption V3, then there exist two sequences of functions
 $\{\,\underline{V}\,_n\}_{n=1}^\infty$ and
 $\,\{\,\overline{V}\,_n\}_{n=1}^\infty$ such that
 $\underline{V}\,_n$ and $\overline{V}\,_n$ satisfy Assumption V$3'$ and the following
 \bee
 &&d\underline{V}\,_{n,\,t}=-\underline{g}\,_{n,\,t}\,dt
 +\underline{Z}\,^k_{n,\,t}\,dB_t^k,\qquad
 d\overline{V}_{n,\,t}=-\overline{g}_{n,\,t}\,dt
 +\overline{Z}\,^k_{n,\,t}\,dB_t^k,
 \\[2mm]
 &&\qquad\quad\underline{V}\,_n\rightarrow \underline{V}\,,\;\;
 \overline{V}\,_n\rightarrow \overline{V}\;\;
 \mbox{a.e. in}\;\;\Omega\times\overline{Q}\;\;
 \mbox{and in}\;\;\mathbb{H}^{1,\,2}\cap\mathbb{S}^{0,\,2}.
 \eee
 Moreover, there exists a sequence of  nonnegative random fields
 $\{\,\widetilde{h}_n\,\}_{n=1}^\infty$ such that
 $\widetilde{h}_n\in \mathbb{H}^{0,\,2}$ and
\be                                                                                                  \label{eq4.2}
 \left.
 \begin{array}{l}
 {\cal L}\,\underline{V}\,_n+{\cal M}^k \underline{Z}\,_n^k-\,
 \underline{g}\,_n
 \geq -\widetilde{h}_n\;,\;\;\;\;
 {\cal L}\,\overline{V}\,_n+{\cal M}^k \overline{Z}\,_n^k
 -\,\overline{g}\,_n
 \leq \widetilde{h}_n\;\;\mbox{a.e. in}\;\;\Omega\times \overline{Q},
 \vspace{4mm} \\
 \|h_n\|_{0,\,2}\leq C(K)\,\Big(\,\|\,\underline{V}\,\|_{1,\,2}
 +\|\,\overline{V}\,\|_{1,\,2}
 +\|\,\underline{Z}\,\|_{0,\,2}+\|\,\overline{Z}\,\|_{0,\,2}
 +\|\,h\,\|_{0,\,2}\,\Big).
 \end{array}
 \right.
\ee
\end{prop}

 \noindent{\bf Proof.} (i) For the two processes $\underline{V}$ and $\overline{V}$ in
 Assumption V$3'$, define
 $$
 h\triangleq\max\Big\{\,\Big(\,{\cal L}\,\underline{V}\,+{\cal M}^k
 \underline{Z}\,^k-\,\underline{g}\,\Big)^-,\,
 \Big(\,{\cal L}\,\overline{V}\,+{\cal M}^k \overline{Z}\,^k-\,\overline{g}\,\Big)^+\,\Big\},
 $$
 where  $f^+$ and $f^-$ represent the positive and negative parts of $f$, respectively.
 %That means that
% $$
% f^+=
% \left\{
% \begin{array}{lll}
% f\;\;&\mbox{ \rm if }\;f>0,
% \vspace{2mm}\\
% 0\;\;&\mbox{ \rm if }\;f\leq0;
% \end{array}
% \right.\qquad
% f^-=
% \left\{
% \begin{array}{lll}
% 0\;\;&\mbox{ \rm if }\;f\geq0,
% \vspace{2mm}\\
% -f\;\;&\mbox{ \rm if }\;f<0.
% \end{array}
% \right.
% $$
 We can check that $\underline{V}$ and $\overline{V}$ satisfy
 Assumption V$3$.\medskip

 (ii) For the process $\underline{V}$ in Assumption V$3$, define
 $$
 \underline{V}\,_n=\,\underline{V}\,*\zeta_n,\qquad
 \underline{g}\,_n=\,\underline{g}\,*\zeta_n,\qquad
 \underline{Z}\,_n=\,\underline{Z}\,*\zeta_n,\qquad
 h_n=\,h*\zeta_n.
 $$
 We have
 $$
 d\underline{V}\,_{n,\,t}=-g\,_{n,\,t}\,dt
 +\underline{Z}\,^k_{n,\,t}\,dB_t^k,\quad
 \underline{V}\,_n\,\in \mathbb{H}^{2,\,2},\,
 \underline{g}\,_n\,\in \mathbb{H}^{0,\,2},\,
 \underline{Z}\,_n\,\in \mathbb{H}^{1,\,2},\,
 h_n\,\in \mathbb{H}^{0,\,2}.
 $$
 Since  $\underline{V}\,\in \mathbb{H}^{1,\,2},\;\underline{g}\,,\,
 \underline{Z}\,,\,\underline{h}\,\in \mathbb{H}^{0,\,2}$ and
 $\underline{V}\,$ has the special representation (\ref{specialform}),
 then we have that
 $$
 \underline{V}\,_n\rightarrow \underline{V}\;\;
 \mbox{a.e. in}\;\;\Omega\times \overline{Q}\;\;
 \mbox{and in}\;\;\mathbb{H}^{1,\,2}\cap\mathbb{S}^{0,\,2},\quad
 \underline{g}\,_n\rightarrow\,\underline{g}\,,\;
 \underline{Z}\,_n\rightarrow\,\underline{Z}\,,\;
 h_n\rightarrow\,h\;\;
 \mbox{in}\;\;\mathbb{H}^{0,\,2}.
 $$
 So, it is sufficient to prove that there exist
 a nonnegative random field sequence  $\{\,\widetilde{h}_n\,\}_{n=1}^\infty$
 satisfying  (\ref{eq4.2}).\smallskip

 Let $\eta\,(x)=\zeta_n(y-x)$ in (\ref{eq4.1}).  In the following, we
 estimate every term in (\ref{eq4.1}). At first, we define
 \bee
 I^1_n(t,y)&\triangleq &-\int_{\mathbb{R}^{d^*}}a^{ij}(t,x)\,
 D_i \,\underline{V}\,(t,x)\,D_j\,\eta(x)\,dx
 -a^{ij}(t,y)D_{i,j}\,\underline{V}\,_n(t,y).
 \eee
 We have
 \bee
 |\,I^1_n(t,y)\,|&\leq& \int_{\mathbb{R}^{d^*}}|\,a^{ij}(t,x)-a^{ij}(t,y)\,|\,
 \,|\,D_i \,\underline{V}\,(t,x)\,n^{d^*+1}\,D_j\,\zeta(\,n(y-x))\,|\,dx
 \\[2mm]
 &\leq& K\,n^{d^*+1}\,\int_{\mathbb{R}^{d^*}}|\,x-y\,|\,
 \,|\,D_i \,\underline{V}\,(t,x)\,D_j\,\zeta(\,n(y-x))\,|\,dx
 \\[2mm]
 &=& K\,\,\int_{B_1}|\,x\,|\,
 \Big|\,D_i \,\underline{V}\,\Big(t,y-{x\over n}\Big)
 \,D_j\,\zeta(x)\,\Big|\,dx
 \\[2mm]
 &\leq& C(K)\,n^{d^*}\,\int_{D_n(y)}
 \,|\,D_i \,\underline{V}\,(t,x)\,|\,dx,
 \eee
 where
 $$
 D_n(y)\stackrel{\Delta}
 {=}\Big\{\,(x_1,\cdot\cdot\cdot,x_{d^*}):|\,x_i-y_i\,|\leq {1\over n}\;\;
 \mbox{for any}\;\; i=1,\cdot\cdot\cdot,d^*\,\Big\}.
 $$
 Hence, we have the following estimate
 \bee
 \|I^1_n\|^2_{0,2}
 &\leq& C(K)\,\mathbb{E}\Bigg[\,\int_Q\,\int_{D_n(y)}
 \,n^{d^*}\,|\,D_i \,\underline{V}\,(t,x)\,|^2\,dx\,dy\,dt\,\Bigg]
 \leq C(K)\,\|D\,\underline{V}\,\|^2_{0,2}\;.
 \eee

 Denote
 \bee
 I^2_n(t,y)&\!\!\!\triangleq \!\!\!&-\int_{\mathbb{R}^{d^*}}
 (\,\sigma^{jk}\,\underline{Z}\,^k\,D_j\,\eta\,)\,(t,x)\,dx
 -(\,\sigma^{jk}\,D_j\,\underline{Z}\,_n^k\,)\,(t,y),
 \\[2mm]
 I^3_n(t,y)&\!\!\!\triangleq \!\!\!&\int_{\mathbb{R}^{d^*}}
 \Big[\Big(b^iD_i\underline{V}+c\underline{V}
 +\mu^k\underline{Z}\,^k \Big)\eta\Big](t,x)dx
 -\Big(b^iD_i\underline{V}\,_n+c\underline{V}\,_n
 +\mu^k\underline{Z}\,_n^k\Big)(t,y),\qquad
 \\[2mm]
 I^4_n(t,y)&\!\!\!\triangleq \!\!\!&-\int_{\mathbb{R}^{d^*}}
 \Big[\,\Big(\,D_j\,a^{ij}\,D_i\,\underline{V}\,
 +D_i\,\sigma^{ik}\,\underline{Z}\,^k\,\Big)\,\eta\,\Big]\,(t,x)\,dx.
 \eee
 Repeating the above argument, we get the following estimate:
 \bee
 \|I^2_n\|^2_{0,2}+\|I^3_n\|^2_{0,2}+\|I^4_n\|^2_{0,2}
 \leq  C(K)\,\Big(\,\|\,\underline{Z}\,\|^2_{0,2}
 +\|\,\underline{V}\,\|^2_{1,2}\,\Big)\;.
 \eee
 Hence, if we denote
 $$
 \widetilde{h}_n=h_n+|\,I^1_n\,|+|\,I^2_n\,|+|\,I^3_n\,|+|\,I^4_n\,|,
 $$
 then the above estimates and
 ${\cal T}(\,\underline{V}\,,\,\underline{Z}\,,\,\underline{G}\,,\,
 h\,,\,\eta\,)\geq0$ imply that
 \bee
 {\cal L}\,\underline{V}\,_n+{\cal M}^k \underline{Z}\,_n^k-\,\underline{g}\,_n
 +\widetilde{h}_n \geq 0,\quad
 \|\widetilde{h}_n\|^2_{0,2}
 \leq  C(K)\,\Big(\,\|\,\underline{Z}\,\|^2_{0,2}
 +\|\,\underline{V}\,\|^2_{1,2}+\|\,h\,\|^2_{0,2}\,\Big)\;.
 \eee

 Using the same method, we can deduce the rest of (\ref{eq4.2}).
 \hfill$\Box$
 \medskip

 A strong solution of BSPDVI (\ref{BSPDI2}) is defined
 as follows.

\begin{defn}                                                                               \label{defn4.1}
 If $(V,Z,k^+,k^-)\in \mathbb{H}^{2,\,2}\times\mathbb{H}^{1,\,2}
 \times \mathbb{H}^{0,\,2}\times\mathbb{H}^{0,\,2}$
 such that
\be                                                                                         \label{eq4.3}
 \left\{
 \begin{array}{l}
 \displaystyle V_t=\varphi
 +\int_t^T\!\!({\cal L} V_s+{\cal M}^k Z^k_s+f_s+k^+_s-k^-_s)\,ds\\[4mm]
 \displaystyle \quad\quad\quad-\int_t^T\!\!\!Z^k_s\,dB^k_s,
 \quad \mbox{\rm a.e. } x\in \mathbb{R}^{d^*}\;
 \mbox{\rm for all } t\in [0,T]\;\mbox{\rm and a.s. in}\;\Omega;
 \vspace{4mm} \\
 \underline{V}\leq V\leq \overline{V}\,,\;\;\;
 k^\pm\geq0\;\;\;\mbox{a.e. in}\;\Omega\times Q;
 \vspace{4mm} \\
 \displaystyle{\int_0^T (V_t-\underline{V}\,_t)\,k_t^+\,dt
 =\int_0^T (\overline{V}_t-V_t)\,k_t^-\,dt=0}
 \;\;\;\mbox{a.e. in}\;\Omega\times\mathbb{R}^{d^*}.
 \end{array}
 \right.
\ee
 Then $(V,Z,k^+,k^-)$ is called a strong solution of BSPDVI (\ref{BSPDI2}).
\end{defn}

 We have the following verification theorem for Problem $\sD_{tx}$.

\begin{thm}                                                                                        \label{thm4.2}
 { \rm ({\bf Verification}) }\, Let Assumptions V1-V4 be satisfied and
 $(t,x)\in Q$. Let $X$ be the solution of
 SDE~\eqref{eq1.1} with the value being $x$ at the initial time $t$, and
 $\underline{V}\,(X)\triangleq \underline{V}\,(\cdot,X_{\cdot})$ and
 $\overline{V}(X)\triangleq \overline{V}\,(\cdot,X_{\cdot})$ are stochastic
 processes with continuous paths. Assume that the four-tuple $(V,Z,k^+,k^-)$
 is a strong solution of BSPDVI~\eqref{BSPDI2}
 with
\be                                                                                                \label{eq4.4}
 a^{ij} \triangleq {1\over2}\left(\sum_{l=1}^{d_1}\gamma_{i\,l}\gamma_{j\,l}
 +\sum_{l=1}^{d_2}\theta_{i\,l}\theta_{j\,l}\right),\quad
 b^i \triangleq \beta_{i},\quad c \triangleq 0,\quad
 \sigma^{ik} \triangleq \theta_{i\,k},\quad
 \mu_k \triangleq 0
\ee
 for any $i,j=1,\cdots, d^*$ and $k=1,\cdots,d_2$.
 Then $V(t,x)$ is the value of Problem $\sD_{tx}$. Define
\bee
 \tau_1^*\triangleq \inf\left\{s\in[\,t,T\,]:V_s(X_s)
 =\underline{V}\,_s(X_s)\right\}\wedge T
 \eee
 and
\bee
 \tau_2^* \triangleq \inf\left\{s\in[\,t,T\,]:V_s(X_s)
 =\overline{V}\,_s(X_s)\right\}\wedge T.
\eee
 Then, $(\tau_1^*, \tau_2^*)$ is a Nash equilibrium point of Problem $\sD_{tx}$.
\end{thm}

\begin{rmk}                                                                                       \label{rem4.1}
 Assumptions D1, D2, and (\ref{eq4.4}) imply
 Assumptions V1 and V2. If  Assumptions D1, D2, V$3'$, and (\ref{eq4.4})
 are all satisfied, then the  processes $\underline{V}\,(X)$ and
 $\overline{V}(X)$ are It\^o processes and possess path continuous
 versions by Theorem~\ref{thm3.1}. If $\underline{V}$ and $\overline{V}$
 are continuous with respect to $(t,x)$ a.s. in $\Omega$, and
 Assumptions D1 and D2 are satisfied, then $\underline{V}\,(X)$
 and $\overline{V}(X)$ are stochastic processes of continuous paths.
\end{rmk}

 \noindent{\bf Proof of Theorem~\ref{thm4.2}.}
 It is sufficient to prove that for any $\tau_1,\,\tau_2\in{\cal U}\,_{t,T}$,
 it holds that
 \bee
 \mathbb{E}\,\Big[\,R_t(x;\tau^*_1,\tau_2)\,\Big|\,{\cal F}_t\,\Big]
 \geq V_t(x)
 \geq \mathbb{E}\,\Big[\,R_t(x;\tau_1,\tau_2^*)\,\Big|\,{\cal F}_t\,\Big]\;\;
 \mbox{a.s. in}\;\;\Omega,
 \eee
 with equality in the first inequality if  $\tau_2=\tau^*_2$ and in the
 second  inequality if $\tau_1=\tau^*_1$. In what follows, we only prove
 the second inequality  since the first one can be proved in a symmetric
 way.\smallskip

 From Theorem~\ref{thm3.1}, we deduce that $V(X)$ is an It\^o process and
 possesses a path continuous version. Hence, $V(X)-\underline{V}\,(X)$ and
 $\overline{V}(X)-V(X)$ are stochastic processes with continuous paths.
 So, we have almost everywhere in $\Omega\times Q$,
 $$
 \chi_{\{\tau_1^*<T\}}\,V_{\tau_1^*}(X_{\tau_1^*})=
 \chi_{\{\tau_1^*<T\}}\,\underline{V}\,_{\tau_1^*}(X_{\tau_1^*}) \mbox{ \rm and }
 \chi_{\{\tau_2^*<T\}}\,V_{\tau_2^*}(X_{\tau_2^*})
 =\chi_{\{\tau_2^*<T\}}\,\overline{V}\,_{\tau_2^*}(X_{\tau_2^*}),
 $$
 and a.s. in $\Omega$,
 $$
 V_s(X_s)-\underline{V}\,_s(X_s)>0\;\mbox{for any}\;t\leq s<\tau_1^*\mbox{ \rm and }
 \overline{V}_s(X_s)-V_s(X_s)>0\;\mbox{for any}\;t\leq s<\tau_2^*.
 $$
 Moreover, the third equality in Definition~\ref{defn4.1} implies that:
 $$
 k^+\,\chi_{\{(w,s,x):\, (V-\underline{V}\,)(w,s,x)>0\,\}}=0,\;\;\;
 k^-\,\chi_{\{(w,s,x):\, (\overline{V}-V)(w,s,x)>0\,\}}=0\;\;\;
 \mbox{in}\;\;\mathbb{M}^{d^*+3}.
 $$
 From Lemma~\ref{lem3.3} , we deduce that
 $$
 k_s^+(X_s)\,\chi_{\{\,s<\tau_1^*\}}=0,\quad
 k_s^-(X_s)\,\chi_{\{\,s<\tau_2^*\}}=0\;\; \mbox{a.e. in}\;\;\Omega\times(0,T),
 $$
 and for any $\tau_1,\,\tau_2\in{\cal U}\,_{t,\,T}$ satisfying
 $\tau_1\leq \tau_1^*,\,\tau_2\leq \tau_2^*$, the following hold
 \be                                                                                            \label{eqadd1}
 \int_t^{\tau_1}\,k_s^+(X_s)\,ds=0,\quad\mbox{and}\quad
 \int_t^{\tau_2}\,k_s^-(X_s)\,ds=0,\;\;\mbox{a.s. in}\;\Omega.
 \ee

 On the event $\{\tau_1\in{\cal U}\,_{t,\,T}:\tau_1\geq \tau_2^*\}$,
 applying Theorem~\ref{thm3.1}, we have
 \bee
 &&R_t(x;\tau_1,\tau_2^*)
 \\[2mm]
 &=&\int_t^{\tau_2^*}f_u(X^{t,x}_u)\,du
 +\overline{V}_{\tau_2^*}(X^{t,x}_{\tau_2^*})\,\chi_{\{\tau_2^*<T\}}
 +\varphi(X^{t,x}_T)\chi_{\{\tau_2^*= T\}}
 \\[2mm]
 &=&\int_t^{\tau_2^*}f_u(X^{t,x}_u)\,du+V_{\tau_2^*}(X^{t,x}_{\tau_2^*})\qquad
 \\[2mm]
 &=&\int_t^{\tau_2^*} f_u(X^{t,x}_u)\,du
 +\int_t^{\tau_2^*}\Bigg[\,L V_u+M^k Z^k_u-({\cal L} V_u
 +{\cal M}^k Z^k_u+f_u+k^+_u-k^-_u)\Bigg]\,(X^{t,x}_u)\,du\qquad
 \\[2mm]
 &&+M_1(\tau_2^*)+V_t(X^{t,x}_t),
 \eee
 where
 $$
 M_1(\tau_2^*)\triangleq
 \int_t^{\tau_2^*} (Z^k_u+M^kV_u)(X^{t,x}_u)\,dB^k_u
 +\int_t^{\tau_2^*} (N^lV_u)\,(X^{t,x}_u)\,dW^l_u.
 $$
 Recalling (\ref{eq4.4}), (\ref{eqadd1}) and $k^+\geq0$, we have
 $$
 R_t(x;\tau_1,\tau_2^*)
 =V_t(x)-\int_t^{\tau_2^*}k^+_u(X^{t,x}_u)\,du+M_1(\tau_2^*)
 \leq V_t(x)+M_1(\tau_2^*) \;\;\mbox{a.s. in}\;\;\Omega,
 $$
 with equality if $\tau_1=\tau^*_1$, which follows from $\tau_2^*\leq\tau_1^*$
 and  (\ref{eqadd1}).\smallskip

 On the event
 $\{\tau\in{\cal U}\,_{t,\,T}:\tau_1<\tau_2^*\}$, in a similar way, we have
 \bee
 R_t(x;\tau_1,\tau_2^*)&=&
 \int_t^{\tau_1} f_s(X^{t,x}_s)\,ds
 +\underline{V}\,_{\tau_1}(X^{t,x}_{\tau_1})
 \leq \int_t^{\tau_1} f_s(X^{t,x}_s)\,ds+V_{\tau_1}(X^{t,x}_{\tau_1})
 \\[2mm]
 &\leq&V_t(x)+M_1(\tau_1)\;\;\mbox{a.s. in}\;\;\Omega,\qquad
 \eee
 with equality if $\tau_1=\tau^*_1$. So, we obtain that
 $$
 R_t(x;\tau_1,\tau_2^*)
 \leq V_t(x)+M_1({\tau_1\wedge \tau_2^*}),
 $$
 almost surely with equality if $\tau_1=\tau^*_1$. Taking conditional
 expectations with respect to ${\cal F}_t$, we have
 $$
 \mathbb{E}\,\Big[\,R_t(x;\tau_1,\tau_2^*)\,\Big|\,{\cal F}_t\,\Big]
 \leq V_t(x)\;\;
 \mbox{for any}\;\tau_1\in{\cal U}\,_{t,\,T}\;\;
 \mbox{a.s. in}\;\;\Omega
 $$
 with the equality holding if $\tau_1=\tau_1^*$. The proof is then
 complete. \hfill$\Box$ \\

\setcounter{equation}{0} \setcounter{section}{4}
\section{Strong solution of BSPDVI
 (\ref{BSPDI2}): existence and uniqueness, and comparison theorem.}

 In this section, we use the penalty method to prove the existence, and establish a comparison
 theorem  to obtain the uniqueness.

 We first deduce some estimates about BSPDE (\ref{eq2.3}), which are
 crucial to the proof of the main results in this paper.

\begin{lem}                                                                                   \label{lem5.1}
 Let the assumptions in Lemma~\ref{lem2.2} be satisfied. Define
 $f(\cdot)\triangleq F(\cdot,V(\cdot))$.

 \noindent Then the  strong solution $(V,Z)$ of BSPDE~(\ref{eq2.3})
 satisfies the following:
\be\nonumber
 &&\|V^-\|^2_{1,\,2}+ |\|V^-|\|^2_{0,\,2}+\|Z\chi_{\{V\leq0\}}\|^2_{0,\,2}
 \\[2mm]                                                                                       \label{eq5.1}
 &\leq& C(\kappa,K,T)\,\mathbb{E}\,\left(\,|\,\varphi^-|^2_{0,\,2}
 +\int_0^T\!\!\!\int_{D}\,f_s^-\,V_s^-\,dx\,ds\right),
 \qquad
 \\[2mm]\nonumber
 &&\|V^+\|^2_{1,\,2}+ |\|V^+|\|^2_{0,\,2}+\|Z\chi_{\{V\geq0\}}\|^2_{0,\,2}
 \\[2mm]                                                                                        \label{eq5.2}
 &\leq& C(\kappa,K,T)\,\mathbb{E}\,\left(\,|\,\varphi^+|^2_{0,\,2}
 +\int_0^T\!\!\!\int_{D}\,f_s^+\,V_s^+\,dx\,ds\right),
 \\[2mm]                                                                                        \label{eq5.3}
 &&\|V\|^2_{1,\,2}+ |\|V|\|^2_{0,\,2}+\|Z\|^2_{0,\,2}
 \leq
 C(\kappa,K,T)\,\mathbb{E}\,\left(\,|\,\varphi|^2_{0,\,2}
 +\int_0^T\!\!\!\int_{D}\,(f_s V_s)^+\,dx\,ds\right),\qquad
 \\[2mm]                                                                                         \label{eq5.4}
 &&\mathbb{E}\,\Bigg(\,\int_0^T\!\!\!\int_D\,f_s\,V_s^-\,dx\,ds\,\Bigg)
 \leq C(\kappa,K)\Big(\,\|V^-\|^2_{0,\,2}+
 \,\mathbb{E}\,|\,\varphi^-\,|_{0,\,2}^2\,\Big),
 \\[2mm]                                                                                        \label{eq5.5}
 &&\mathbb{E}\,\Bigg(\,\int_0^T\!\!\!\int_D\,\,-f_s\,V_s^+\,dx\,ds\,\Bigg)
 \leq C(\kappa,K)\Big(\,\|V^+\|^2_{0,\,2}
 +\,\mathbb{E}\,|\,\varphi^+\,|_{0,\,2}^2\,\Big).
\ee
\end{lem}

 \noindent{\bf Proof.} Define the auxiliary functions:
 $$
 \xi(r)=
 \left\{
 \begin{array}{ll}
 r^2,&r\leq 0;
 \vspace{2mm} \\
 r^2\,(1-r)^3,\quad&0\leq r\leq1;
 \vspace{2mm} \\
 0,\quad&r\geq1
 \end{array}
 \right.
 $$
 and
 $$
 \xi_n(r)={1\over n^2}\;\xi(nr).
 $$
we have $\xi_n\in C\,^2(\mathbb{R}),\;|\xi''_n|\leq C$,
 $$
  \lim\limits_{n\rightarrow\infty}\xi_n(r)=(r^-)^2\;\mbox{ \rm uniformly},\quad
 \lim\limits_{n\rightarrow\infty}\xi'_n(r)=-2r^-\;\mbox{ \rm uniformly},
 $$
and
$$
 \lim\limits_{n\rightarrow\infty}\xi''_n(r)
 =\left\{
 \begin{array}{ll}
 2,\quad &r\leq0;
 \vspace{2mm} \\
 0, &r>0.
 \end{array}
 \right.
 $$

 Applying It\^o's formula for Hilbert-valued semimartingales (see e.g. in
 \cite{Prato}) to $\xi_n(V)$, we deduce that $J^1_n=J^2_n-J^3_n\;$
 a.s. in $\Omega$, where
\bee
 J^1_n&\triangleq &\int_D\,\xi_n(V_t)\,dx
 -\int_D\,\xi_n(\varphi)\,dx\,,
 \\[2mm]
 J^2_n&\triangleq &\int_D\,\int_t^T\,
 \Bigg[\,\xi'_n( V_s)\,  D_i(a^{ij} D_jV_s+\sigma^{ik} Z^k_s)
 +\xi'_n(V_s)\,\Big(\,(b^i- D_ja^{ij})\,  D_i V_s
 \\[2mm]
 &&+c\, V_s +(\mu^k-  D_i\sigma^{ik})\, Z^k_s
 + f_s\,\Big)-{1\over2}\,\xi''_n(V_s)\,|Z_s|^2
 \,\Bigg]\,ds\,dx
\eee \bee
 &=&\int_D\,\int_t^T\,
 \Bigg[\,-\,\xi''_n( V_s)\,\Big(\,a^{ij}  D_iV_s D_j V_s
 +\sigma^{ik} Z^k_s  D_i V_s+{1\over2}\,|Z_s|^2\Big)
 \\[2mm]
 &&+\xi'_n( V_s)\,\Big(\,(b^i- D_ja^{ij})\,  D_iV_s
 +c\,V_s+(\mu^k-  D_i\sigma^{ik})\,Z^k_s
 +f_s\,\Big)\,\Bigg]\,ds\,dx,
 \\[2mm]
 J^3_n&\triangleq &\int_D\,\int_t^T\,\xi'_n(V_s)\,
 Z_s^k\,d\,B_s^k\,dx.
\eee

 Taking $n\rightarrow \infty$ in the above equalities, we have that
\bee
 J^1_n&\rightarrow&J_1\triangleq \int_D\,(\,V_t^-\,)\,^2\,dx
 -\int_D\,(\,\varphi^-\,)\,^2\,dx,
 \\[2mm]
 J^2_n&\rightarrow&J_2\triangleq \int_D\,\int_t^T\,
 \chi_{\{ V_s\leq0\}}\,\Bigg[\,-\Big(\,2a^{ij}  D_i V_s D_j V_s
 +2\sigma^{ik} Z^k_s  D_i V_s+|Z_s|^2\,\Big)
 \\[2mm]
 &&-2\,V_s^-\,\Big(\,(b^i- D_ja^{ij})\,  D_i V_s
 +c\, V_s+(\mu^k-  D_i\sigma^{ik})\,Z^k_s
 + f_s\,\Big)\,\Bigg]\,ds\,dx
 \\[2mm]
 J^3_n&\rightarrow&J_3\triangleq
 -\,2\,\int_D\,\int_t^T\,V_s^-
 \,Z_s^k\,d\,B_s^k\,dx.
\eee
 Hence, $J_1=J_2-J_3\;$ a.s. in $\Omega$. Denote
 $$
 |\|V|\|^2_{k,\,2\,;\,t}\triangleq
 E\left(\sup\limits_{s\in[\,t,\,T\,]}|\,V_s\,|_{k,\,2}^2\right)\quad \hbox{ \rm and }\quad
 \|V\|^2_{k,\,2\,;\,t}\triangleq
 E\left(\,\int_t^T\,|\,V_s\,|_{k,\,2}^2\,ds\right).
 $$
 Since
 \bee
 J_2\leq\int_D\,\int_t^T\!\!\!\!\!
 \chi_{\{ V_s\leq0\}}\Bigg[\!\!-{\kappa\over2}\,|D V_s|^2
 +C(\kappa,K)\,|\,V_s|^2-{\kappa\over 2\,(\kappa+4K)}\,|Z_s|^2
 -2\,V_s^-\,f_s\,\Bigg]ds\,dx,\qquad
 \eee
 taking expectation on both sides of $J_1=J_2-J_3$, we have
 \be\nonumber
 \mathbb{E}|V_t^-|_{0,\,2}^2-\mathbb{E}|\varphi^-|_{0,\,2}^2
 &\leq&-{\kappa\over2}\,\|V^-\|^2_{1,\,2\,;\,t}
 -{\kappa\over 2\,(\kappa+4K)}\,\|Z\chi_{\{V\leq0\}}\|^2_{0,\,2\,;\,t}
 \\[2mm]                                                                                    \label{eq5.6}
 &&+\,C(\kappa,K)\,\|V^-\|^2_{0,\,2\,;\,t}
 -2 \mathbb{E}\,\left(\int_t^T\!\int_D\,V_s^-\,f_s\,dx\,ds\right).
 \ee
 Hence,
 \bee
 &&\sup\limits_{s\in[\,t,T\,]}\mathbb{E}|V_s^-|_{0,\,2}^2
 +\|V^-\|^2_{1,\,2\,;\,t}+\|Z\chi_{\{V\leq0\}}\|^2_{0,\,2\,;\,t}
 \\[2mm]
 &\leq& C(\kappa,K)\,\left\{\,\mathbb{E}\,(\,|\,\varphi^-\,|_{0,\,2}^2\,)
 +\|V^-\|^2_{0,\,2\,;\,t}+
 \mathbb{E}\,\left(\int_t^T\int_D\,f_s^-\,V_s^-\,dx\,ds\right)\right\}.
 \eee
 This along with the Gronwall inequality yields that
 \bee
 &&\sup\limits_{t\in[\,0,T\,]}\mathbb{E}|V_t^-|_{0,\,2}^2
 +\|V^-\|^2_{1,\,2}+\|Z\chi_{\{V\leq0\}}\|^2_{0,\,2}
 \\[2mm]
 &\leq& C(\kappa,K,T)\;\mathbb{E}\,\left(\,|\,\varphi^-\,|_{0,\,2}^2\,
 +\int_0^T\int_D\,f_s^-\,V_s^-\,dx\,ds\right).
 \eee
 Using the BDG inequality, we have
 \bee
 &&|\|V^-|\|^2_{0,\,2}+\|V^-\|^2_{1,\,2}+\|Z\chi_{\{V\leq0\}}\|^2_{0,\,2}
 \\[2mm]
 &\leq&{1\over2}\,|\|V^-|\|^2_{0,\,2}
 +C(\kappa,K,T)\left[\,\|Z\chi_{\{V\leq0\}}\|^2_{0,\,2}+
 \mathbb{E}\,\left(|\,\varphi^-\,|_{0,\,2}^2
 +\int_0^T\,\int_D\,f_s^-\,V_s^-\,dx\,ds\right)\,\right].\qquad
 \eee
 So, the Gronwall inequality implies (\ref{eq5.1}).

 We prove (\ref{eq5.2}) in the same way. Since
 $$
 |\,V\,|+|\,DV\,|=V^++V^-+|\,DV^+\,|+|\,DV^-\,|,\quad
 (f\,V)^+=f^+\,V^++f^-\,V^-,
 $$
 we derive (\ref{eq5.3}) from (\ref{eq5.1}) and (\ref{eq5.2}).

 On the other hand, taking $t=0$ in (\ref{eq5.6}), we
 have (\ref{eq5.4}). Applying (\ref{eq5.4}) to the BSPDE for $(-V,-Z)$, we have (\ref{eq5.5}).
 \hfill$\Box$

 We have the following comparison theorem, which implies the uniqueness of strong solution of BSPDVI (\ref{BSPDI2}).

\begin{thm}                                                                                 \label{thm5.2}
 Let Assumptions V1 and V2 be satisfied.
 Let $(V_i,Z_i,k_i^+,k_i^-)$ be the strong solution to BSPDVI~\eqref{BSPDI2}
  with  $(f, \,\varphi, \,\overline{V},\,\underline{V})\triangleq (f_i,\,\varphi_i,\,\overline{V}_i,\,\underline{V}\,_i)$ for $i=1,2$.
 If $f_1\geq f_2,\,\varphi_1\geq \varphi_2,\, \overline{V}_1\geq \overline{V}_2,$
 and $ \underline{V}\,_1\geq \underline{V}\,_2$, then $V_1\geq V_2$ a.e.
 in $\Omega\times Q$.
\end{thm}

\begin{rmk}  \label{rem5.1} In general, we have no comparison on the reflection
parts $k^\pm$. Here is an illustration.  Take $d^*=1, {\cal L} =D_x^2,
 {\cal M}=D_x,$ and $ T=1$.

(i) For the following parameters:
\bee
 &&\overline{V}\,_1=(1+x^2)^{-1},\qquad \underline{V}\,_1=0,\qquad
 \varphi_1=(1+x^2)^{-1},\qquad f_1=3\,(1+x^2)^{-3};
  \\[2mm]
 &&\overline{V}\,_2=(1+x^2)^{-1},\qquad \underline{V}\,_2=0,\qquad
 \varphi_2=0,\qquad f_2=-(1+x^2)^{-3},\eee
 we have
 \bee
 && V_1=(1+x^2)^{-1},\qquad Z_1=0,\qquad  k_1^+=0,\qquad
 k_1^-=(1+6x^2)\,(1+x^2)^{-3};
 \\[2mm]
 && V_2=0,\qquad Z_2=0,\qquad  k_2^+=(1+x^2)^{-3},\qquad
 k_2^-=0,
\eee
 with
 $$  k_1^+<k_2^+\quad\hbox{ \rm and } \quad  k_1^->k_2^-.$$

 (ii) For the following parameters:
\bee
 &\!\!\!\!\!\!&\overline{V}\,_1=e^4\,(1+x^2)^{-1},\quad
 \underline{V}\,_1=e^{3-3t}\,(1+x^2)^{-1},\quad
 \varphi_1=e^3\,(1+x^2)^{-1},\quad f_1=0;\qquad
 \\[2mm]
 &\!\!\!\!\!\!&\overline{V}\,_2=e^{3t-3}\,(1+x^2)^{-1},\quad
 \underline{V}\,_2=0,\qquad
 \varphi_2=e^{-3}\,(1+x^2)^{-1},\qquad f_2=0,\qquad
\eee
 we have
\bee
 &\!\!\!\!\!\!& V_1=e^{3-3t}\,(1+x^2)^{-1},\quad Z_1=0,\quad
 k_1^+=e^{3-3t}\,(3x^4+5)\,(1+x^2)^{-3},\quad
 k_1^-=0;
 \\[2mm]
 &\!\!\!\!\!\!& V_2=e^{3t-3}\,(1+x^2)^{-1},\;\; Z_2=0,\;\;
 k_2^+=0,\;\;
 k_2^-=e^{3-3t}\,(3x^4+12x^2+1)\,(1+x^2)^{-3},\qquad
\eee with
 $$k_1^+>k_2^+\quad\hbox{ \rm and } \quad  k_1^-<k_2^-.$$
\end{rmk}

 \noindent{\bf Proof of Theorem~\ref{thm5.2}.} Define
 \bee
 &&\Delta V \triangleq  V_1-V_2,\qquad\;\;\;
 \Delta Z \triangleq  Z_1-Z_2,\qquad
 \Delta k^+ \triangleq  k^+_1-k_2^+,\qquad\qquad\qquad\\[2mm]
 &&\Delta k^- \triangleq  k^-_1-k_2^-,\qquad
 \Delta f \triangleq  f_1-f_2,\qquad\;\;
 \Delta \varphi \triangleq  \varphi_1-\varphi_2.
 \eee
 Then $\Delta f\geq0,\;\Delta \varphi\geq0$ and $\Delta V$ satisfies the
 following BSPDE:
 $$
 \Delta V_t=\Delta \varphi+\int_t^T\,({\cal L}\Delta V_s+{\cal M}^k\Delta Z^k_s+\Delta f_s
 +\Delta k_s^+-\Delta k_s^-)\,ds-\int_t^T\,\Delta Z^k_s\,dB^k_s.
 $$

 In view of Lemma~\ref{lem5.1} , since $\Delta f\geq0,\;\Delta \varphi\geq0,\,
 k^\pm_1\geq0,\,k^\pm_2\geq0$, we have that
 \be\nonumber
 &&\|(\Delta V)^-\|^2_{1,\,2}+
 |\|(\Delta V)^-|\|^2_{0,\,2}+\|Z\chi_{\{\Delta V\leq0\}}\|^2_{0,\,2}
 \\[2mm]\nonumber
 &\leq&C\,\mathbb{E}\,\left(\,|\,(\Delta \varphi)^-|^2_{0,\,2}
 +\int_0^T\,\int_{\mathbb{R}^{d^*}}\,
 (\Delta f_s+\Delta k_s^+-\Delta k_s^-)^-(\Delta V_s)^-\,dx\,ds\right)
 \\[2mm]                                                                                        \label{eq5.7}
 &\leq& C\,\mathbb{E}\,\left(\,
 \int_{\mathbb{R}^{d^*}}\,\int_0^T\,\chi_{\{V_1<V_2\}}\,
 (\,k_{1,\,s}^-+k_{2,\,s}^+\,)\,(V_{2,\,s}-V_{1,\,s})\,ds\,dx\right).
 \ee

 On the other hand, in the domain $\{V_1<V_2\}$, we have
 $$
 \underline{V}\,_2\leq \underline{V}\,_1\leq V_1<V_2\leq
 \overline{V}\,_2\leq \overline{V}\,_1.
 $$
 Hence, we deduce that
 \be\nonumber
 &&\int_0^T\,\chi_{\{V_1<V_2\}}\,(\,k_{1,\,s}^-+k_{2,\,s}^+\,)
 \,(V_{2,\,s}-V_{1,\,s})\,ds
 \\[2mm]\nonumber
 &=&\int_0^T\,\chi_{\{V_1<V_2\}}\,
 \Big[\,(\,V_{2,\,s}-\overline{V}\,_{1,\,s}\,)\,k^-_{1,\,s}
 -(\,V_{1,\,s}-\overline{V}\,_{1,\,s}\,)\,k^-_{1,\,s}\,\Big]\,ds
 \\[2mm]\nonumber
 &&+\int_0^T\,\chi_{\{V_1<V_2\}}\,
 \Big[\,(\,V_{2,\,s}-\underline{V}\,_{2,\,s}\,)\,k^+_{2,\,s}
 -(\,V_{1,\,s}-\underline{V}\,_{2,\,s}\,)\,k^+_{2,\,s}\,\Big]\,ds
 \\[2mm]                                                                                     \label{eq5.8}
 &\leq&\int_0^T\,\chi_{\{V_1<V_2\}}\,
 \Big[\,-(\,V_{1,\,s}-\overline{V}\,_{1,\,s}\,)\,k^-_{1,\,s}
 +(\,V_{2,\,s}-\underline{V}\,_{2,\,s}\,)\,k^+_{2,\,s}\,\Big]\,ds.
 \ee

 Since $(V_1-\overline{V}_1)\,k^-_1\leq0$ and
 $(V_2-\underline{V}\,_2)\,k^+_2\geq0$ a.e. in $\Omega\times Q$, then the
 forth equality in (\ref{eq4.3}) implies that
 \be                                                                                         \label{eq5.9}
 \int_0^T\,\chi_{\{V_1<V_2\}}\,(V_{1,\,s}-\overline{V}\,_{1,\,s})\,k^-_{1,\,s}\,ds
 =\int_0^T\,\chi_{\{V_1<V_2\}}\,(\,V_{2,\,s}-\underline{V}\,_{2,\,s}\,)\,k^+_{2,\,s}\,ds
 =0.
 \ee

 From (\ref{eq5.7})-(\ref{eq5.9}), we conclude that
 $$
 \|(\Delta V)^-\|^2_{1,\,2}+
 |\|(\Delta V)^-|\|^2_{0,\,2}+\|Z\chi_{\{\Delta V\leq0\}}\|^2_{0,\,2}=0,
 $$
 which means that $\Delta V\geq0$, i.e., $V_1\geq V_2$.
 \hfill$\Box$

The main result of this section is stated as follows.

\begin{thm}                                                                                        \label{thm5.3}
 Let Assumptions V1-V4 be satisfied. Then BSPDVI~\eqref{BSPDI2} has a unique
 strong solution $(V,Z,k^+,k^-)$ such that
\be\nonumber
 &\!\!\!\!\!\!&\|V\|_{2,\,2}+ |\|V|\|_{1,\,2}+\|Z\|_{1,\,2}
 +\|k^+\|_{0,\,2}+\|k^-\|_{0,\,2}
 \\[2mm]                                                                                           \label{eq5.10}
 &\!\!\!\leq\!\!\!& C(\kappa,K,T)\,\Big(\,\mathbb{E}\,[\,|\,\varphi\,|\,_{1,\,2}\,]
 +\|f\|_{0,\,2}+\|\,\underline{V}\,\|_{1,\,2}+\|\,\overline{V}\,\|_{1,\,2}
 \nonumber\\
 &&\quad\qquad +\|\,\underline{Z}\,\|_{0,\,2}+\|\,\overline{Z}\,\|_{0,\,2}
 +\|\,h\,\|_{0,\,2}\,\Big).\qquad
\ee
 Moreover, if $\underline{V}\,(X)$ and $\overline{V}(X)$ are stochastic
 processes of continuous paths, then the strong solution of  BSPDVI~\eqref{BSPDI2}
 coincides with the value of Problem $\sD_{tx}$.
\end{thm}

 \noindent{\bf Proof.}  The uniqueness of strong solutions is an immediate
 consequence of Theorem~\ref{thm5.2}. The rest of the proof is divided into the four steps.

 {\bf Step 1. Penalty. }  We construct the following penalized problem to
 approximate BSPDE~(\ref{BSPDI2}):
\be \label{eq5.11}
 \left\{
 \begin{array}{l}
 dV_{n,\,t}=-\left[{\cal L} V_{n,\,t}+{\cal M}^k Z^k_{n,\,t}+f_t
 +n\,(\underline{V}\,_{n,\,t}-V_{n,\,t})^+
 -n\,(V_{n,\,t}-\overline{V}_{n,\,t}\,)^+\right]\,dt\\[4mm]
 \quad\quad \quad\quad +Z^k_{n,\,t}\,dB^k_t;
 \vspace{3mm} \\
 V_{n,\,T}(x)=\varphi_n(x)\,,
 \end{array}
 \right.
\ee
 where $\underline{V}\,_n$ and $\overline{V}_n$ are defined in Proposition~\ref{prop4.1},
 and $\varphi_n$ is the mollification of $\varphi$, i.e., $\varphi_n=\varphi*\zeta$,
 which is defined in (\ref{eq3.1}). It is clear that
 $\underline{V}\,_n\leq\varphi_n\leq\overline{V}_n$ and $\varphi_n$ converges
 to $\varphi$ in $\mathbb{L}^{1,2}$.\smallskip

 Denote
 $$
 F(w,t,x,u)
 \stackrel{\triangle}{=}
 f_t(w,x)+n\,(\,\underline{V}\,_{n,\,t}(w,x)-u\,)^+
 -n\,(\,u-\overline{V}_{n,\,t}(w,x)\,)^+.
 $$
 Then $F(\cdot,u)$ is
 ${\cal P}^B\times{\cal B}(\mathbb{R}^{d^*})$-measurable for any
 $u\in \mathbb{R}$ and
 $$
 F(\cdot,0)=f+\,n\,\underline{V}\,^+_n -n\,\overline{V}\,^-_n\in \mathbb{H}^{0,\,2},
 \quad |F(w,t,x,u)-F(w,t,x,v)|\leq 2n|\,u-v\,|.
 $$
 According to Lemma~\ref{lem2.2}, BSPDE~(\ref{eq5.11}) has a unique
 strong solution denoted by $(V_n, Z_n)$.\smallskip

 {\bf Step 2. The sequence $\{V_n\}$ is bounded in  $\mathbb{H}^{2,\,2}$. } Define
 $$
 \Delta\underline{V}\,_n \stackrel{\triangle}{=}
 \sqrt{n}\,(\,V_n-\,\underline{V}\,_n),\qquad
 \Delta\underline{Z}\,_n \stackrel{\triangle}{=}
 \sqrt{n}\,(\,Z_n-\,\underline{ Z}\,_n),\qquad
 \Delta\overline{ V}\,_n \stackrel{\triangle}{=}
 \sqrt{n}\,(\,V_n-\,\overline{ V}\,_n).
 $$
 Then $(\Delta\underline{V}\,_n, \Delta\underline{Z}\,_n)$ is a strong solution to the following BSPDE:
 \bee
 \left\{
 \begin{array}{ll}
 d\,\Delta\underline{V}\,_{n,\,t}=-({\cal L}\,\Delta \underline{V}\,_{n,\,t}
 +{\cal M}^k \Delta \underline{Z}\,^k_{n,\,t}+\sqrt{n}\,\underline{f}\,_{n,\,t}
 +n\,(\Delta \underline{ V}\,_{n,\,t})^-
 -n\,(\Delta \overline{ V}\,_{n,\,t})^+\,)\,dt
 \\[2mm]
 \hspace{10.1cm}+\;\Delta \underline{Z}\,\,^k_{n,\,t}\,dB^k_t;\\
 \Delta\,\underline{V}\,_{n,\,T}(x)=\underline{\varphi}\,_n(x)
 \end{array}
 \right.
 \eee
 with
 $$
 \underline{f}\,_n\stackrel{\triangle}{=}f+{\cal L}\,\underline{V}\,_n
 +{\cal M}^k \underline{Z}\,_n^k-\underline{g}\,_n\in \mathbb{H}^{0,\,2} \quad \hbox{\rm and } \quad
 \underline{\varphi}\,_n\stackrel{\triangle}{=}
 \sqrt{n}\,(\,\varphi_n-\,\underline{V}\,_{n,\,T}\,)\geq0.
 $$

 From Lemma~\ref{lem2.5}, we have $\underline{\varphi}\,_n\in \mathbb{H}^{1,\,2}$.
 In view of  (\ref{eq4.2}), we have
 $$
 \underline{f}\,_n\geq f-\,\widetilde{h}_n,\qquad
 \underline{f}\,^-_n\leq f^-+\,\widetilde{h}_n,
 $$
 $$
 \|\,\underline{f}\,^-_n\|_{0,\,2}\leq
 C(K)\,\Big(\,\|\,f\,\|_{0,\,2}+\|\,\underline{V}\,\|_{1,\,2}
 +\|\,\overline{V}\,\|_{1,\,2}+\|\,\underline{Z}\,\|_{0,\,2}
 +\|\,\overline{Z}\,\|_{0,\,2}+\|\,h\,\|_{0,\,2}\,\Big).
 $$
Since
 $$
 \Big(\,\sqrt{n}\,\underline{f}\,_n+n\,(\Delta \underline{ V}\,_n)^-
 -n\,(\Delta \overline{ V}\,_n)^+\,\Big)^-
 \leq \sqrt{n}\,\underline{f}\,_n^-
 +n\,(\Delta \overline{ V}\,_n)^+,\;\;
 (\Delta \overline{ V}\,_n)^+\,(\Delta \underline{ V}\,_n)^-=0,
 $$
 we have from (\ref{eq5.1}) in Lemma~\ref{lem5.1}  that
 \be\nonumber
 &&\|(\,\Delta\underline{V}\,_n)^-\|^2_{1,\,2}
 +|\|(\,\Delta\underline{V}\,_n)^-|\|^2_{0,\,2}
 +\|\Delta \underline{Z}\,_n\chi_{\{\Delta\underline{V}\,_n\leq0\}}\|^2_{0,\,2}
 \\[2mm]\nonumber
 &\leq&C(\kappa,K,T)\,\mathbb{E}\,\left[\,
 \int_0^T\,\int_{\mathbb{R}^{d^*}}\,
 \Big(\,\sqrt{n}\,\underline{f}\,_{n,\,s}^-+n\,(\Delta \overline{ V}\,_{n,\,s})^+\,\Big)
 \,(\,\Delta\underline{V}\,_{n,\,s})^-\,dx\,ds\right]
 \\[2mm]\nonumber
 &=&C(\kappa,K,T)\,\mathbb{E}\,\left[\,
 \int_0^T\,\int_{\mathbb{R}^{d^*}}\,
 \underline{f}\,_{n,\,s}^-\,
 (\,\sqrt{n}\,\Delta\underline{V}\,_{n,\,s})^-\,dx\,ds\right]
 \\[2mm]                                                                                            \label{eq5.12}
 &\leq&{n\over 4\,C(\kappa,K)}\;\|
 \,(\,\Delta\underline{V}\,_n)^-\,\|^2_{0,\,2}
 +C(\kappa,K,T)\,C(\kappa,K)\,\|\,\underline{f}\,_n^-\,\|^2_{0,\,2}\;.
 \ee

 Recalling  $(\Delta \overline{ V}\,_n)^+\,(\Delta \underline{ V}\,_n)^-=0$,
 and applying (\ref{eq5.4}) in Lemma~\ref{lem5.1} to BSPDE of
 $(\Delta\underline{V}\,_n, \Delta\underline{Z}\,_n)$, we have that
 \bee
 &&\mathbb{E}\,\Bigg[\,\int_0^T\,\int_{\mathbb{R}^{d^*}}\,
 \sqrt{n}\,\underline{f}\,_{n,\,s}\,
 (\,\Delta \underline{ V}\,_{n,\,s}\,)^-
 +n\,\Big((\Delta \underline{V}\,_{n,\,s})^-\,\Big)^2\,dx\,ds\,\Bigg]
 \\[2mm]
 &=& \mathbb{E}\,\Bigg[\,\int_0^T\,\int_{\mathbb{R}^{d^*}}\,
 \Big(\,\sqrt{n}\,\underline{f}\,_{n,\,s}
 +n\,(\Delta \underline{V}\,_{n,\,s})^-
 -n\,(\Delta \overline{V}\,_{n,\,s})^+\,\Big)\,
 (\,\Delta \underline{V}\,_{n,\,s}\,)^-\,dx\,ds\,\Bigg]
 \\[2mm]
 &\leq&  C(\kappa,K)\,\Big(\,\|(\,\Delta\underline{V}\,_n)^-\|^2_{0,\,2}
 +\mathbb{E}\,[\,|\,\underline{\varphi}\,_n^-\,|_{0,\,2}^2\,]\,\Big)
 \leq C(\kappa,K)\,\|(\,\Delta\underline{V}\,_n)^-\|^2_{1,\,2}.
 \eee
 Hence, we deduce that
 \bee
 &&n^2\;\|\,(\,V_n-\,\underline{V}\,_n\,)^-\,\|^2_{0,\,2}
 =\mathbb{E}\,\Bigg[\,\int_0^T\,\int_{\mathbb{R}^{d^*}}\,
 n\,\Big((\Delta \underline{V}\,_{n,\,s})^-\,\Big)^2\,dx\,ds\,\Bigg]
 \\[2mm]
 &\leq&  C(\kappa,K)\,\|(\,\Delta\underline{V}\,_n)^-\|^2_{1,\,2}
 -\mathbb{E}\,\Bigg[\,\int_0^T\,\int_{\mathbb{R}^{d^*}}
 \underline{f}\,_{n,\,s}\,(\sqrt{n}\,
 \Delta \underline{V}\,_{n,\,s})^-\,dx\,ds\,\Bigg]\qquad
 \\[2mm]
 &\leq&{n^2\over 2}\;\|\,(\,V_n-\,\underline{V}\,_n\,)^-\,\|^2_{0,\,2}
 +C(\kappa,K,T)\,\|\,\underline{f}\,_n^-\,\|^2_{0,\,2}\,.\qquad
 \qquad\qquad\qquad(\mbox{by}\;(\ref{eq5.12}))
 \eee
 Hence, we have
  $$
 \|\,k_n^+\,\|^2_{0,\,2}
 \leq C(\kappa,K,T)\,\|\,\underline{f}\,_n^-\,\|^2_{0,\,2}
 $$
 with
 $$
 k_n^+\stackrel{\triangle}{=}n\,(\,\underline{V}\,_n-V_{n}\,)^+.
 $$

 In a similar way, we have
 \bee
 &\!\!\!\!\!\!&\|\,k_n^-\,\|^2_{0,\,2}
 \leq C(\kappa,K,T)\,\|\,f+{\cal L}\,\overline{V}\,_n
 +{\cal M}^k \overline{Z}\,_n^k-\overline{g}\,_n\,\|^2_{0,\,2}
 \\[2mm]
 &\!\!\!\leq\!\!\!& C(\kappa,K,T)\,\Big(\,\|\,\underline{V}\,\|^2_{1,\,2}
 +\|\,\overline{V}\,\|^2_{1,\,2}
 +\|\,\underline{Z}\,\|^2_{0,\,2}+\|\,\overline{Z}\,\|^2_{0,\,2}
 +\|\,h\,\|^2_{0,\,2}+\|\,f\,\|^2_{0,\,2}\,\Big)
 \eee
 with
 $$k_n^-\stackrel{\triangle}{=}n\,(\,V_{n}-\overline{V}\,_n\,)^+.
 $$

 From Lemma~\ref{lem2.2}, we have the following estimate:
 \be\nonumber
 &\!\!\!\!\!\!&\|V_{n}\|^2_{2,\,2}+ |\|V_{n}|\|^2_{1,\,2}
 +\|Z_{n}\|^2_{1,\,2}
 +\|k_n^+\|^2_{0,\,2}+\|k_n^-\|^2_{0,\,2}
 \\[2mm]\nonumber
 &\!\!\!\leq\!\!\!& C(\kappa,K,T)\,
 \Big(\mathbb{E}\,[\,|\,\varphi_n\,|\,^2_{1,\,2}\,]
 +\|f+k_{n}^+-k_{n}^-\|^2_{0,\,2}
 +\|k_n^+\|^2_{0,\,2}+\|k_n^-\|^2_{0,\,2}\Big)
 \\[2mm]\nonumber
 &\!\!\!\leq\!\!\!& C(\kappa,K,T)\,
 \Big(\,\mathbb{E}\,[\,|\,\varphi\,|\,^2_{1,\,2}\,]
 +\|f\|^2_{0,\,2}+\|\,\underline{V}\,\|^2_{1,\,2}
 +\|\,\overline{V}\,\|^2_{1,\,2}
 \\[2mm]                                                                                               \label{eq5.13}
 &&\qquad\qquad\quad +\|\,\underline{Z}\,\|^2_{0,\,2}
 +\|\,\overline{Z}\,\|^2_{0,\,2}
 +\|\,h\,\|^2_{0,\,2}\,\Big).
 \ee
 Hence, there exists a subsequence of $\{(V_{n},Z_{n},k^+_n,\,k^-_n)\}$,
 still denoted by itself, such that $V_{n},\,Z_{n},\,k^+_n,\,k^-_n$
 converges to $V,Z,k^+,$ and $k^-$  weakly in $\mathbb{H}^{2,\,2},\,\mathbb{H}^{1,\,2}
 ,\,\mathbb{H}^{0,\,2},$ and $\mathbb{H}^{0,\,2}$, respectively.
Letting $n$ tend to $\infty$ in (\ref{eq5.13}), we have (\ref{eq5.10}).
 \smallskip

{\bf Step 3. $V_{n}$ converges to $V$ strongly in
 $\mathbb{H}^{1,\,2}$. }  It is sufficient to show that
 $\{V_{n}\}_{n=1}^{\infty}$ is a Cauchy sequence in $\mathbb{H}^{1,\,2}$. Define
 $$
 \Delta V_{n,\,m}\triangleq  V_n-V_m,\;\;
 \Delta \underline{V}\,_{n,\,m}\triangleq \underline{V}\,_n-\underline{V}\,_m,\;\;
 \Delta \overline{V}\,_{n,\,m}\triangleq \overline{V}\,_n-\overline{V}\,_m,\;\;
 \Delta Z_{n,\,m} \triangleq  Z_n-Z_m,
 $$
 and $\Delta k_{n,\,m}^\pm \triangleq  k^\pm_n-k^\pm_m$ with $m\geq n$. Then
 $(\Delta V_{n,\,m}, \Delta Z_{n,\,m})$ satisfies the following BSPDE:
 \bee
 \Delta V_{n,\,m,\,t}=\int_t^T({\cal L} \Delta V_{n,\,m,\,s}+{\cal M}^k \Delta Z^k_{n,\,m,\,s}
 +\Delta k_{n,\,m,\,s}^+-\Delta k_{n,\,m,\,s}^-\,)\,ds
 -\int_t^T\Delta Z^k_{n,\,m,\,s}\,dB^k_s.
 \eee

 From (\ref{eq5.3}) in Lemma~\ref{lem5.1}, we have
 \be  \nonumber
 &&\|\,\Delta V_{n,\,m}\,\|_{1,\,2}^2+|\|\,\Delta V_{n,\,m}\,|\|_{0,\,2}^2
 +\|\,\Delta Z_{n,\,m}\,\|^2_{0,\,2}
 \\[2mm]                                                                                          \label{eq5.14}
 &\leq&C\,\mathbb{E}\,\Bigg\{\,\int_0^T\,\int_{\mathbb{R}^{d^*}}\,
 \Big[\,\Delta V_{n,\,m,\,s}\,(\,\Delta k_{n,\,m,\,s}^+
 -\Delta k_{n,\,m,\,s}^-\,) \,\Big]^+\,dx\,ds\,\Bigg\}.
 \ee
 Moreover, we have
 \bee
 \Delta V_{n,\,m}\,\Delta k_{n,\,m}^+
 &=&n\,[\,(\,V_n-\underline{V}\,_n\,)-(\,V_m-\underline{V}\,_m\,)
 +\Delta \underline{V}\;_{n,\,m}\,]\,(\,\underline{V}\,_n-V_n\,)^+
 \\[2mm]
 &&- m\,[\,(\,V_n-\underline{V}\,_n\,)-(\,V_m-\underline{V}\,_m\,)
 +\Delta \underline{V}\;_{n,\,m}\,]\,(\,\underline{V}\,_m-V_m\,)^+
 \\[2mm]
 &\leq&(\,\underline{V}\,_m-V_m\,)^+\,k^+_n+
 (\,\underline{V}\,_n-V_n\,)^+\,k^+_m
 +\Delta \underline{V}\;_{n,\,m}\,(\,k^+_n-\,k^+_m\,)
 \\[2mm]
 &\leq& {2\over n}\;k^+_n\,k^+_m+
 |\,\Delta \underline{V}\;_{n,\,m}\,|\,(\,|\,k^+_n\,|+|\,\,k^+_m\,|\,)\eee
 and
 \bee
 \Delta V_{n,\,m}\,\Delta k_{n,\,m}^-
 &=&n\,[\,(\,V_n-\overline{V}\,_n\,)-(\,V_m-\overline{V}\,_m)
 +\Delta \overline{V}\;_{n,\,m}\,]\,(\,V_n-\overline{V}\,_n\,)^+
 \\[2mm]
 &&-m\,[\,(\,V_n-\overline{V}\,_n\,)-(\,V_m-\overline{V}\,_m)
 +\Delta \overline{V}\;_{n,\,m}\,]\,(\,V_m-\overline{V}\,_m\,)^+\qquad
 \\[2mm]
 &\geq&-\,(\,V_m-\overline{V}\,_m\,)^+\,k^-_n
 -(\,V_n-\overline{V}\,_n\,)^+\,k^-_m
 +\Delta\overline{V}\,_{n,\,m}\,(\,k^-_n-\,k^-_m\,)
 \\[2mm]
 &\geq& -{2\over n}\,k^-_n\,k^-_m-
 |\,\Delta \overline{V}\,_{n,\,m}\,|\,(\,|\,k^-_n\,|+|\,\,k^-_m\,|\,).
 \eee
 Combining (\ref{eq5.13}) and (\ref{eq5.14}), we derive that as
 $n,\,m\rightarrow\infty$,
 $$
 \|\Delta V_{n,\,m}\|_{1,\,2}^2+|\|\Delta V_{n,\,m}|\|_{0,\,2}^2
 + \|\Delta Z_{n,\,m}\|^2_{0,\,2}
 \leq {C\over n}
 +C\,(\|\,\Delta \underline{V}\,_{n,\,m}\,\|^2_{0,\,2}
 +\|\,\Delta \overline{V}\,_{n,\,m}\,\|^2_{0,\,2})\rightarrow 0,
 $$
 since $\{\underline{V}\,_n\}_{n=1}^\infty$ and
 $\{\overline{V}\,_n\}_{n=1}^\infty$ are Cauchy sequences in $\mathbb{H}^{0,2}$.
 Hence,  $V_{n}$ converges to $V$ in $\mathbb{H}^{1,2}$.\smallskip

 Since
 $$
 \|\,(\,\underline{V}\,_n-V_{n}\,)^+\,\|_{0,\,2}
 ={1\over n}\;\|\,k_n^+\,\|_{0,\,2} \leq C/n\rightarrow 0\;\;
 \mbox{as}\;\;n\rightarrow\infty,
 $$
 we derive that $(\,\underline{V}\,_n-V_{n}\,)^+$ converges to $0$ in
 $\mathbb{H}^{0,2}$. On the other hand, $\underline{V}\,_n-V_{n}$
 converges to $\underline{V}\,-V$ in $\mathbb{H}^{0,2}$. Hence,
 we deduce that $(\,\underline{V}\,-V\,)^+=0$, i.e.,
 $V\geq \underline{V}\,$ a.e. in $\Omega\times Q$. In a similar way,
 we have
 $$
 \underline{V}\leq V \leq \overline{V},\;\;\; k^+,\,k^-\geq0\;\;\;
 \mbox{\rm a.e. in }\;\Omega\times Q,
 $$
 which is the second and third inequalities in (\ref{eq4.3}).\smallskip

 {\bf Step 4. It remains to  show that $(V,\,Z,\,k^+,\,k^-)$ satisfies the
 first and forth equalities in (\ref{eq4.3}). }  Let $\xi$ be an arbitrary
 ${\cal P}^B$-predictable element of  $\mathcal{L}^2$, and
 $\forall\;\eta\in H^{2,\,2}$. We have $\eta\,\xi\in \mathbb{H}^{2,\,2}$.
 \smallskip

 Rewrite (\ref{eq5.11}) into the integrated form:
 $$
 V_{n,\,t}=\varphi+\int_t^T({\cal L} V_{n,\,s}+{\cal M}^k Z^k_{n,\,s}+f_s
 +k_{n,\,s}^+-k_{n,\,s}^-\,)\,ds -\int_t^T Z^k_{n,\,s}\,dB^k_s.
 $$
 Multiplying  by $\eta\,\xi$ both sides of the last equality, and integrating, we have
 $$
 \mathbb{E}\left(\,\int_0^T
 \int_{\mathbb{R}^{d^*}}V_{n\,,t}(x)\,\eta(x)\,\xi_t\,dx\,dt\,\right)
 =I_{1}+I_{2},
 $$
 with
 \bee
 I_{1}&\triangleq &\mathbb{E}\left(\,\int_0^T \int_{\mathbb{R}^{d^*}}
 \varphi(x)\,\eta(x)\,\xi_t\,dx\,dt\,\right)
 -\mathbb{E}\left(\,\int_0^T \int_{\mathbb{R}^{d^*}}
 \int_t^T \eta(x)\,\xi_t\,Z^k_{n,\,s}(x)\,dB^k_s\,dx\,dt\,\right),\qquad
 \\[2mm]
 I_{2}&\triangleq &\mathbb{E}\left[\,\int_0^T \int_{\mathbb{R}^{d^*}}\int_t^T\eta(x)\,\xi_t\,
 \Big({\cal L} V_{n,\,s}+{\cal M}^k Z^k_{n,\,s}+f_s+k_{n,\,s}^+-k_{n,\,s}^-\,\Big)\,(x)
 \,ds\,dx\,dt\,\right].
 \eee

 Evidently, we have
 $$
\lim_{n\to \infty} \mathbb{E}\left(\,\int_0^T \int_{\mathbb{R}^{d^*}}
 V_{n,\,t}(x)\,\eta(x)\,\xi_t\,dx\,dt\,\right)
 = \mathbb{E}\left(\,\int_0^T\,\xi_t \int_{\mathbb{R}^{d^*}}
 V_t(x)\,\eta(x)\,dx\,dt\,\right).
 $$
 Since $Z_{n}$ converges to $Z$  weakly in $\mathbb{H}^{1,\,2}$, from a
 known result (see \cite[Theorem 4, page 63]{Rozovskii}), we have
 $$
 \lim_{n\to \infty}\mathbb{E}\left(\,\xi_t\,\int_t^T\,\int_{\mathbb{R}^{d^*}}\eta(x)
 \,Z^k_{n,\,s}(x)\,dx\,dB^k_s\,\right)
 = \mathbb{E}\left(\,\xi_t\,\int_t^T\,\int_{\mathbb{R}^{d^*}}
 \eta(x)\,Z^k_s(x)\,dx\,dB^k_s\,\right)
 $$
 for any $t\in[\,0,T\,].$  Moreover,  for every $t\in[\,0,T\,]$,
 \bee
 &&\Bigg|\,\mathbb{E}\left(\,\xi_t\,\int_t^T\,\int_{\mathbb{R}^{d^*}}\eta(x)
 \,Z^k_{n,\,s}(x)\,dx\,dB^k_s\,\right)\,\Bigg|\\[2mm]
 &=&\int_{\mathbb{R}^{d^*}}|\,\eta(x)\,|\;\mathbb{E}\,\Bigg|\,\xi_t\,\int_t^T\,
 \,Z^k_{n,\,s}(x)\,dB^k_s\,\Bigg|\,dx
 \\[2mm]
 &\leq&|\,\eta\,|\,_{0,2}\,
 \left\{\,\int_{\mathbb{R}^{d^*}}\,\mathbb{E}\,(\,\xi^2_t\,)\;
 \mathbb{E}\,\left[\;\left|\;\int_t^T Z^k_{n,\,s}(x)\,
 dB^k_s\,\right|^2\;\right]\,dx\;\right\}^{1\over2}
 \\[2mm]
 &\leq& |\,\eta\,|\,_{0,2}\,
 \left\{\,\int_{\mathbb{R}^{d^*}}\,\mathbb{E}\,(\,\xi^2_t\,)\;
 \mathbb{E}\,\left[\;\int_t^T \left|\;Z^k_{n,\,s}(x)\,\right|^2\,
 ds\;\right]\,dx\;\right\}^{1\over2}\\
 &\leq& |\,\eta\,|\,_{0,2}\Big[\,\mathbb{E}\,(\,\xi^2_t\,)\,\Big]^{1/2}
 \|Z_{n}\|^2_{0,\,2}
 \leq C\,|\,\eta\,|\,_{0,2}\Big[\,\mathbb{E}\,(\,\xi^2_t\,)\,\Big]^{1/2}.
 \eee
 Hence, using Lebesgue's dominant convergence theorem, we have
 \bee
 \mathbb{E}\left(\,\int_0^T \int_{\mathbb{R}^{d^*}}\int_t^T
 \eta(x)\,\xi_t\,Z^k_{n,\,s}(x)\,dB^k_s\,dx\,dt\,\right)
 \rightarrow\mathbb{E}\left(\,\int_0^T\,\xi_t\int_t^T\, \int_{\mathbb{R}^{d^*}}
 \eta(x)Z^k_{s}(x)\,dx\,dB^k_s\,dt\,\right).
 \eee
In a similar way, we have
 $$
 I_{2}\rightarrow
 \mathbb{E}\left\{\,\int_0^T \xi_t\,\int_t^T\int_{\mathbb{R}^{d^*}}
 \eta(x)\,\Big[\,{\cal L} V_{s}(x)+{\cal M}^k Z^k_{s}(x)+f_s(x)+k_s^+(x)-k_s^-(x)\,\Big]\,
 dx\,ds\,dt\,\right\}.
 $$
 Hence, we show that for any $\eta\in H^{2,\,2}$ and ${\cal P}^B$-predictable
 stochastic process $\xi$ belonging to $\mathcal{L}^2$, it holds that
 \bee
 &&\mathbb{E}\left[\,\int_0^T \xi_t\,\int_{\mathbb{R}^{d^*}}V_t(x)\,\eta(x)\,dx\,dt\,\right]
 \\[2mm]
 &=&\mathbb{E}\left[\,\int_0^T \xi_t\,\left(\int_{\mathbb{R}^{d^*}}\varphi(x)\,\eta(x)\,dx
 -\int_t^T\,\int_{\mathbb{R}^{d^*}}\eta(x)\,Z^k_{s}(x)\,dx\,dB^k_s
 \,\right)\,dt\,\right]
 \\[2mm]
 &+&\mathbb{E}\left[\,\int_0^T \xi_t\int_t^T\int_{\mathbb{R}^{d^*}}\eta(x)\,
 ({\cal L} V_{s}(x)+{\cal M}^k Z^k_{s}(x)+f_s(x)+k_s^+(x)-k_s^-(x)\,)\,dx\,ds\,dt\,\right].
 \eee

 Since $\xi$ is arbitrary, we see that for any $\eta\in H^{2,\,2}$,
 a.e. in $\Omega\times[\,0,T\,]$, it holds that
 $$
 \int_{\mathbb{R}^{d^*}}\!\!\eta\,V_t\,dx
 = \int_{\mathbb{R}^{d^*}}\!\!\eta\,\varphi\,dx
 +\int_t^T\int_{\mathbb{R}^{d^*}}\!\!\eta\,({\cal L} v_{s}
 +{\cal M}^k Z^k_{s}+f_s+k_s^+-k_s^-\,)\,dx\,ds
 -\int_t^T\,\int_{\mathbb{R}^{d^*}}\!\!\eta\,Z^k_{s}\,dx\,dB^k_s.
 $$
 From Lemma~\ref{lem2.5}, we deduce that $V\in \mathbb{H}^{2,\,2}\cap
 \mathbb{S}^{1,\,2}$, and $(V,Z,k^+,k^-)$ satisfies the first equality
 in (\ref{eq4.3}).
 \smallskip

 Now, we prove the forth equality in (\ref{eq4.3}). We have
 $$
 V\geq\,\underline{V}\,,\quad
 k^+\geq0,\quad
 \int_0^T (V_t-\underline{V}\,_t)\,k_t^+\,dt\geq0.
 $$
 On the other hand, since $k^+_n$  converges to $k^+$ weakly in
 $\mathbb{H}^{0,2},\, V_n-\underline{V}\,_n$  converges to
 $V-\underline{V}\,$ strongly in $\mathbb{H}^{0,2}$, and
 $$
 \int_0^T (V_{n,\,t}-\underline{V}\,_{n,\,t})\,k^+_{n,\,t}\,dt
 =n\,\int_0^T (V_{n,\,t}-\underline{V}\,_{n,\,t})\,
 (\underline{V}\,_{n,\,t}-V_{n,\,t})^+\,dt\leq0.
 $$
 Therefore, we have
 $$
 \mathbb{E}\left[\,\int_{\mathbb{R}^{d^*}}\!\!
 \int_0^T (V_t-\underline{V}\,_t)\,k_t^+\,dt\,dx\,\right]
 \leq0,
 $$
 which implies that
 $$
 \int_0^T (V_t-\underline{V}\,_t)\,k_t^+\,dt=0\;\;
 \mbox{a.e. in}\;\Omega\times \mathbb{R}^{d^*}.
 $$
 In a similar way, we have
 $$
 \int_0^T (\overline{V}\,_t-V_t)\,k_t^-\,dt=0\;\;
 \mbox{a.e. in}\;\Omega\times \mathbb{R}^{d^*}.
 $$
 \hfill$\Box$

\setcounter{equation}{0} \setcounter{section}{5}
\section{Further properties on the strong solution, and the free boundary of
 BSPDVI (\ref{BSPDI2}).}

 In this section, we use the comparison theorem
 for BSPDVI to derive properties of strong solutions and
 define the stochastic free boundary of BSPDVI~\eqref{BSPDI2} with extra conditions.
 For this purpose, consider the following further assumptions.\medskip

 \noindent{\bf Assumption V5(i).} The coefficient functions $a, b, c, \sigma,$
 and $\mu$ are independent of the variable $x_i$ for
 $i\in\{\,1,2,\cdot\cdot\cdot,d^*\,\}$. Moreover, the dimension of the
 state space $d^*<4$. \smallskip

 \noindent{\bf Assumption V6(i).} The functions $\underline{f}\,,\,\Delta{V},$
 and $\,\underline{\varphi}\,$ are increasing (resp. decreasing) in $x_i$, where
 $$
 \underline{f}\,\triangleq f+{\cal L}\,\underline{V}\,
 +{\cal M}^k\,\underline{Z}\,^k-\underline{g}\,,\quad
 \Delta{V}\,\triangleq \,\overline{V}\,-\underline{V}\,,\quad
 \underline{\varphi}\,\triangleq \varphi-\underline{V}\,_T,\quad
 i\in\{\,1,2,\cdot\cdot\cdot,d^*\,\}.
 $$

 \noindent{\bf Assumption V7(i).} The functions $\overline{f}\,,\,
 -\Delta{V},$ and $\,\overline{\varphi}\,$ are increasing (resp. decreasing)
 in $x_i$, where
 $$
 \overline{f}\,\triangleq f+{\cal L}\,\overline{V}\,
 +{\cal M}^k\,\overline{Z}\,^k-\overline{g}\,,\quad
 \overline{\varphi}\,\triangleq \varphi-\overline{V}\,_T,\quad
 i\in\{\,1,2,\cdot\cdot\cdot,d^*\,\}.
 $$

We have

\begin{thm}                                                                                    \label{thm6.1}
 Let Assumptions V1, V2, V3$'$, V4, V5(i),
 and V6(i) be satisfied. Let $(V,Z,k^+,k^-)$ be the unique strong solution of
 BSPDVI~(\ref{BSPDI2}).

 (i)$\;$Then $V-\,\underline{V}\,$ is continuous and
 increasing  (resp. decreasing) with respect to $x_i$  for any
 $x_{\!\not{\, i}}\in\mathbb{R}^{d^*-1}$, a.e. in $\Omega\times[\,0,T\,]$.

(ii)$\;$Define the lower free boundary as
 $$
 \underline{S}\,_i(w,t,x_{\!\not{\, i}})
 \triangleq \sup\,\{x_i:(V-\,\underline{V}\,)(w,t,x)=0\}
 $$
 $$
 ({\rm resp.}\;\; \underline{S}\,_i(w,t,x_{\!\not{\, i}})
 \triangleq  \inf\,\{x_i:(V-\,\underline{V}\,)(w,t,x)=0\})
 $$
 with $x_{\!\not{\, i}} \triangleq  (x_1,\cdot\cdot\cdot,x_{i-1},x_{i+1},
 \cdot\cdot\cdot,x_{d^*})\in \mathbb{R}^{d^*-1},$ and the convention that
 $\sup\O=-\infty$ and $\inf\O=\infty$. Then, we have
 \bee
 V>\,\underline{V}\;\;\mbox{ \rm a.e. in}\;\;
 \{x_i>\underline{S}\,_i(w,t,x_{\!\not{\, i}})\}\quad \hbox{ \rm and }\quad
 V=\,\underline{V}\;\;\mbox{ \rm a.e. in}\;\;
 \{x_i\leq\underline{S}\,_i(w,t,x_{\!\not{\, i}})\}.
 \\[2mm]
 ({\rm resp.}\;\;V>\,\underline{V}\;\;\mbox{ \rm a.e. in}\;\;
 \{x_i<\underline{S}\,_i(w,t,x_{\!\not{\, i}})\}\quad\hbox{ \rm and }\quad
 V=\,\underline{V}\;\;\mbox{ \rm a.e. in}\;\;
 \{x_i\geq\underline{S}\,_i(w,t,x_{\!\not{\, i}})\}.)
 \eee

 (iii) If Assumptions V5(j) and V6(j)
 with $j\neq i$ are further satisfied,  the free boundary $\underline{S}\,_i$
 is monotone in $x_j$  for any
 $x_{\!\not{\, i,}\not{\, j}}\in\mathbb{R}^{d^*-2}$, a.e. in $\Omega\times[\,0,T\,]$.
\end{thm}

\begin{rmk}                                                                                        \label{rem6.1}
 The lower free boundary is the interface between $\{V=\underline{V}\,\}$
 and $\{\underline{V}\,<V<\overline{V}\,\}$, and the upper free boundary in
 Theorem \ref{thm6.2} is the interface between $\{V=\overline{V}\,\}$ and
 $\{\underline{V}\,<V<\overline{V}\,\}$.
\end{rmk}

 \noindent{\bf Proof.} Denote $\delta V \triangleq  V-\,\underline{V}$ and
 $\delta Z \triangleq  Z-\,\underline{Z}\,$.  Then
 $(\delta V,\,\delta Z,\,k^+,\,k^-)$ is the strong solution of the following BSPDVI:
 \bee
 \left\{
 \begin{array}{l}
 d\delta V_t=-({\cal L} \delta V_t+{\cal M}^k\delta Z^k_t+\underline{f}\,_t)\,dt
 +\delta Z^k_t\,dB^k_t,\qquad
 \mbox{if}\;\;0<\delta V_t<\Delta V_t\,;
 \vspace{2mm} \\
 d\delta V_t\leq-({\cal L} \delta V_t+{\cal M}^k \delta Z^k_t+\underline{f}\,_t)\,dt
 +\delta Z^k_t\,dB^k_t,\qquad
 \mbox{if}\;\;\delta V_t=0\,;
 \vspace{2mm} \\
 d\delta V_t\geq-({\cal L} \delta V_t+{\cal M}^k\delta Z^k_t+\underline{f}\,_t)\,dt
 +\delta Z^k_t\,dB^k_t,\qquad
 \mbox{if}\;\;\delta V_t=\Delta V_t\,;
 \vspace{2mm} \\
 \delta V_T(x)=\underline{\varphi}\,(x)\,.
 \end{array}
 \right.
 \eee

 For any fixed $\ep>0$, denote
 \bee
 \widetilde{\delta V}(w,t,x)\triangleq  \delta V(w,t,x+\ep e_i),\;\;
 \widetilde{\delta Z}(w,t,x)\triangleq  \delta Z(w,t,x+\ep e_i),
 \\[2mm]
 \widetilde{k^\pm}(w,t,x) \triangleq  k^\pm(w,t,x+\ep e_i),\quad
 \widetilde{\underline{f}}\,(w,t,x)\triangleq \underline{f}\,(w,t,x+\ep e_i),
 \\[2mm]
 \widetilde{\Delta V}(w,t,x)\triangleq \Delta V(w,t,x+\ep e_i),\;\;
 \widetilde{\underline{\varphi}}\,(w,t,x)\triangleq
 \underline{\varphi}\,(w,t,x+\ep e_i),
 \eee
 where $e_i\triangleq  (0,\cdot\cdot\cdot,0,1,0\cdot\cdot\cdot,0)$ is the $i$-th
 standard coordinate vector. So, Assumption V5(i) implies that
 $(\widetilde{\delta V},\,\widetilde{\delta Z},\,\widetilde{k^+},\,\widetilde{k^-})$
 is the strong solution of the following BSPDVI:
 \bee
 \left\{
 \begin{array}{l}
 d\widetilde{\delta V}_t=-({\cal L} \widetilde{\delta V}_t
 +{\cal M}^k\widetilde{\delta Z}^k_t+\widetilde{\underline{f}}\,_t)\,dt
 +\widetilde{\delta Z}^k_t\,dB^k_t,\qquad
 \mbox{if}\;\;0<\widetilde{\delta V}_t<\widetilde{\Delta V}\,_t\,;
 \vspace{2mm} \\
 d\widetilde{\delta V}_t\leq-({\cal L} \widetilde{\delta V}_t
 +{\cal M}^k\widetilde{\delta Z}^k_t+\widetilde{\underline{f}}\,_t)\,dt
 +\widetilde{\delta Z}^k_t\,dB^k_t,\qquad
 \mbox{if}\;\;\widetilde{\delta V}_t=0\,;
 \vspace{2mm} \\
 d\widetilde{\delta V}_t\geq-({\cal L} \widetilde{\delta V}_t
 +{\cal M}^k\widetilde{\delta Z}^k_t+\widetilde{\underline{f}}\,_t)\,dt
 +\widetilde{\delta Z}^k_t\,dB^k_t,\qquad
 \mbox{if}\;\;\widetilde{\delta V}_t=\widetilde{\Delta V}\,_t\,;
 \vspace{2mm} \\
 \widetilde{\delta V}_T(x)=\widetilde{\underline{\varphi}}\,(x)\,.
 \end{array}
 \right.
 \eee
 Moreover, Assumption V6(i) implies that
 $$
 \widetilde{\underline{f}}\,\geq\underline{f}\,,\;\;
 \widetilde{\Delta V}\,\geq\Delta V\,,\;\;
 \widetilde{\underline{\varphi}}\,\geq\underline{\varphi}\,.\;\;
 ({\rm resp.}\;\;\widetilde{\underline{f}}\,\leq\underline{f}\,,\;\;
 \widetilde{\Delta V}\,\leq\Delta V\,,\;\;
 \widetilde{\underline{\varphi}}\,\leq\underline{\varphi}\,.)
 $$
 In view of Theorem~\ref{thm5.2}, we deduce that $\widetilde{\delta V}\geq\,
 (\,{\rm resp.}\,\leq\,)\,\delta V$
 a.e. in $\Omega\times[\,0,T\,]\times \mathbb{R}^{d^*}$,
 which means that for any $\ep>0,\,(V-\underline{V}\,)(w,t,x+\ep e_i)
 \geq\,(\,{\rm resp.}\,\leq\,)\,(V-\underline{V}\,)(w,t,x)$ a.e.
 in $\Omega\times[\,0,T\,]\times \mathbb{R}^{d^*}$.\smallskip

 Since $V,\,\underline{V}\,\in \mathbb{H}^{2,\,2}$ and $d^*<4$, then the Sobolev
 imbedding theorem implies that $V$ and $\underline{V}\,$ have continuous versions
 such that they are continuous with respect to $x$ a.e. in
 $\Omega\times[\,0,T\,]$. Hence, Conclusion (i) has been proved.\smallskip

 Since $V-\,\underline{V}\,\geq0$, Conclusion (ii) is clear. In the following,
 we prove the last result.\smallskip

 Without loss of generality, we suppose $d=2,\,i=1,\,j=2$ and
 $\underline{f}\,,\,\Delta{V},$ and $\underline{\varphi}$ are increasing
 with respect to $x_1$ and $x_2$. Moreover, we fix $(w,t)\in\Omega\times(0,T)$.
 Then it is sufficient to prove that the free boundary $\underline{S}\,_1(w,t,x_2)$
 is decreasing in $x_2$.
 \smallskip

 According to Conclusion (ii),  for any fixed
 $x_2^0\in \mathbb{R}$,
 $$
 V(w,t,x_1,x_2^0)-\underline{V}\,(w,t,x_1,x_2^0)>0,\quad \forall\;
 x_1>\underline{S}\,_1(w,t,x_2^0).
 $$
 Moreover, since $V-\underline{V}\,$ is increasing in $x_2$, we
 have
 $$
 V(w,t,x_1,x_2)-\underline{V}\,(w,t,x_1,x_2)>0,\quad \forall\;
 x_1>\underline{S}\,_1(w,t,x_2^0),\;x_2>x_2^0.
 $$
 In view of the definition of $\underline{S}\,_1$, we have that
 for any $x_2>x_2^0$,
 $$
 \underline{S}\,_1(w,t,x_2)
 \triangleq \sup\,\{x_1:(V-\,\underline{V}\,)(w,t,x_1,x_2)=0\}
 \leq \underline{S}\,_1(w,t,x_2^0).
 $$
 Then the proof is complete.
 \hfill$\Box$ \medskip

 In a similar way, we have

\begin{thm}                                                                                       \label{thm6.2}
 Let Assumptions V1, V2, V3$'$, V4, V5(i), and V7(i) be satisfied. Let
 $(V,Z,k^+,k^-)$ be the unique strong solution of BSPDVI~(\ref{BSPDI2}).

 (i)$\;$Then $V-\,\overline{V}\,$ is continuous and increasing
 (resp. decreasing) with respect to $x_i$  a.e. in
 $\Omega\times[\,0,T\,]\times \mathbb{R}^{d^*-1}$.

(ii)$\;$Define the upper free boundary
 as
 $$
 \overline{S}\,_i(w,t,x_{\!\not{\, i}})
 \triangleq \inf\,\{x_i:(V-\,\overline{V}\,)(w,t,x)=0\}.
 $$
 $$
 ({\rm resp.}\;\;\overline{S}\,_i(w,t,x_{\!\not{\, i}})
 \triangleq  \sup\,\{x_i:(V-\,\overline{V}\,)(w,t,x)=0\}.)
 $$
 Then we have
 \bee
 V=\,\overline{V}\;\;\mbox{ \rm a.e. in}\;\;
 \{x_i\geq\overline{S}\,_i(w,t,x_{\!\not{\, i}})\}\quad \hbox{ \rm and } \quad
 V<\,\overline{V}\;\;\mbox{ \rm a.e. in}\;\;
 \{x_i<\overline{S}\,_i(w,t,x_{\!\not{\, i}})\}.
 \\[2mm]
 ({\rm resp.}\;\;V=\,\overline{V}\;\;\mbox{ \rm a.e. in}\;\;
 \{x_i\leq\overline{S}\,_i(w,t,x_{\!\not{\, i}})\}\quad\hbox{ \rm and } \quad
 V<\,\overline{V}\;\;\mbox{ \rm a.e. in}\;\;
 \{x_i>\overline{S}\,_i(w,t,x_{\!\not{\, i}})\}.)
 \eee

 (iii) If Assumptions V5(j) and V7(j)
 with $j\neq i$ are further satisfied, the free boundary $\overline{S}\,_i$
 is monotone in $x_j$  for any
 $x_{\!\not{\, i,}\not{\, j}}\in\mathbb{R}^{d^*-2}$, a.e. in $\Omega\times[\,0,T\,]$.
\end{thm}

\setcounter{equation}{0} \setcounter{section}{6}
\section{The optimal stopping time problem as an extreme case of a Dynkin's game.}

 In this section, we consider an optimal stopping time problem  (denoted by Problem $\mathscr{O}$ hereafter), which  involves only one choice variable of stopping
 times.  We show that Problem $\mathscr{O}$
 is a special case of Dynkin games under suitable conditions and identify
 the corresponding results about Problem $\mathscr{O}$
 and BSPDVI with one obstacle.

 The state $X$ is governed by SDE (\ref{eq1.1}).
 The payoff is defined by
 \bee
 P_t\,(x;\tau)=\int_t^{\tau}f_u(X^{t,x}_u)\,du
 +\underline{V}\,_{\tau}(X^{t,x}_{\tau})\,\chi_{\{\tau<T\}}
 +\varphi(X^{t,x}_T)\,\chi_{\{\tau\geq T\}},\quad \tau\in {\cal U}\,_{t,T}.
 \eee
 The optimal stopping problem $\mathscr{O}_{tx}$, associated to the initial data $(t,x)$,
 is to find a stopping time
 $\tau^*\in {\cal U}\,_{t,T}$ such that
 \bee
 \mathbb{E}\,\Big[\,P_t(x;\tau^*)\,\Big|\,{\cal F}_t\,\Big]
 =V_t(x)\triangleq 
 \mathop{{\rm ess.sup}}_{\tau\in{\cal U}\,_{t,T}}
 \mathbb{E}\,\Big[\,P_t(x;\tau)\,\Big|\,{\cal F}_t\,\Big].
 \eee
 The random variable $V(t,x)$ is called the value of Problem $\mathscr{O}_{tx}$.

 Consider the following two assumptions on  the cost functions $f,\,\underline{V}\,,$ and
 $\varphi$.

 \noindent{\bf Assumption O1.} (Regularity)
 $\;f\in \mathbb{H}^{0,\,2},\, \varphi\in \mathbb{L}^{1,\,2}$ and
 $\,\underline{V}\,$ is in the form of
 $$
 d\underline{V}\,_t=-\underline{g}\,_t\,dt+\underline{Z}\,^k_t\,dB_t^k,
 $$
 where $\underline{V}\,\in \mathbb{H}^{2,\,2},\,
 \underline{Z}\,\in \mathbb{H}^{1,\,2},$ and $
 \underline{g}\,\in \mathbb{H}^{0,\,2}$.
 \medskip

 \noindent{\bf Assumption  O2.} (Compatibility)
 $\;\underline{V}\,_T\leq \varphi$.
 \medskip

The HJB equation for  Problem $\mathscr{O}_{tx}$ is  the following BSPDVI with one obstacle:
 \begin{equation}\label{BSPDVI1}
 \left\{
 \begin{array}{l}
 dV_t=-({\cal L} V_t+{\cal M}^k Z^k_t+f_t)\,dt+Z^k_t\,dB^k_t,\qquad
 \mbox{if}\;\;V_t>\underline{V}\,_t\,;
 \vspace{2mm} \\
 dV_t\leq-({\cal L} V_t+{\cal M}^k Z^k_t+f_t)\,dt+Z^k_t\,dB^k_t,\qquad
 \mbox{if}\;\;V_t=\underline{V}\,_t\,;
 \vspace{2mm} \\
 V_T(x)=\varphi(x)\,,
 \end{array}
 \right.
 \end{equation}
 where the operators ${\cal L}$ and ${\cal M}$ are defined by (\ref{eq1.3}).

\begin{defn}                                                                                       \label{defn7.1}
 A triplet  $(V,Z,k^+)\in \mathbb{H}^{2,\,2}\times\mathbb{H}^{1,\,2}\times
 \mathbb{H}^{0,\,2}$   is called a strong solution of BSPDVI~\eqref{BSPDVI1} if it satisfies the following:
 \bee
 \left\{
 \begin{array}{l}
 \displaystyle{V_t=\varphi
 +\int_t^T\!\!({\cal L} V_s+{\cal M}^k Z^k_s+f_s+k^+_s)\,ds
 -\int_t^T\!\!\!Z^k_s\,dB^k_s}\;\;\mbox{a.e. in}\;\mathbb{R}^{d^*}\;
 \mbox{for all }t\;\mbox{a.s. in}\;\Omega;
 \vspace{4mm} \\
 V\geq \underline{V}\,,\;\;\;
 k^+\geq0\;\;\;\mbox{a.e. in}\;\Omega\times \overline{Q};
 \vspace{4mm} \\
 \displaystyle{\int_0^T (V_t-\underline{V}\,_t)\,k_t^+\,dt=0}
 \;\;\;\mbox{a.e. in}\;\Omega\times\mathbb{R}^{d^*}.
 \end{array}
 \right.
 \eee
\end{defn}

 Identical to the proof of Theorem~\ref{thm5.2}, we have the following comparison theorem.

\begin{thm}     \label{thm7.1}
 Let Assumptions V1 and V2 be satisfied. Let
 $(V_i,Z_i,k_i^+)$ be the strong solution of BSPDVI~\eqref{BSPDVI1}
 associated with $(f_i,\,\varphi_i,\,\underline{V}\,_i)$ for $i=1,2$. If
 $f_1\geq f_2,\,\varphi_1\geq \varphi_2,$ and $\underline{V}\,_1\geq \underline{V}\,_2$,
 then $V_1\geq V_2$ a.e. in $\Omega\times Q$.
\end{thm}

 The following lemma gives the relationship between Problems $\sD_{tx}$ and $\mathscr{O}_{tx}$,
 and between BSPDVIs~\eqref{BSPDI2} and~\eqref{BSPDVI1}.

\begin{lem}                                                                                          \label{lem7.2}
 Let Assumptions D1 and D2 (resp. V1 and V2), O1 and O2 be satisfied.
Then there exists a stochastic fields $\overline{V}$ such that Assumptions
 V3$'$ and V4 are satisfied. Moreover, Problems $\mathscr{O}_{tx}$ and $\sD_{tx}$ (resp. BSPDVIs~\eqref{BSPDVI1} and~\eqref{BSPDI2}) are equivalent. And we have the following
 estimate
\be\nonumber
 &&\|\,\overline{V}\,\|_{2,\,2}+\|\,\overline{Z}\,\|_{1,\,2}
 +\|\,\overline{g}\,\|_{0,\,2}
 \\[2mm]                                                                                                \label{eq7.1}
 &\leq& C(\kappa,K,T)\,\Big(\,\mathbb{E}\,[\,|\,\varphi\,|\,_{1,\,2}\,]
 +\|f\|_{0,\,2}+\|\,\underline{V}\,\|_{2,\,2}
 +\|\,\underline{Z}\,\|_{1,\,2}+\|\,\underline{g}\,\|_{0,\,2}+1\,\Big).\quad
\ee
\end{lem}

 \noindent{\bf Proof.} Let Assumptions D1, D2, O1, and O2 be satisfied.
 Let $(\widetilde{V},\,\widetilde{Z})$ be the strong solution of the
 following BSPDE:
\be                                                                                                   \label{eq7.2}
 \left\{
 \begin{array}{l}
 d\widetilde{V}_t=-(L\,\widetilde{V}_t
 + M^k\,\widetilde{Z}^k_t+\widetilde{f}_t)\,dt
 + \widetilde{Z}^k_t\,dB^k_t;
 \vspace{2mm} \\
 \widetilde{V}_T=\varphi^+\,,
 \end{array}
 \right.
\ee
 where $L$ and $M$ are defined by (\ref{eq2.2}) and
 $$
 \widetilde{f}\triangleq\max\{\,f,\,0,\,g-L\,\underline{V}\,
 -M^k\,\underline{Z}\,^k\,\}\in \mathbb{H}^{0,\,2}.
 $$

 According to Theorem 2.2 in \cite{Du}, BSPDE (\ref{eq7.2}) has a strong solution
 $(\widetilde{V},\,\widetilde{Z})\in \mathbb{H}^{2,\,2}\times\mathbb{H}^{1,\,2}$.
 Moreover, the comparison theorem for linear BSPDE in \cite{Du} implies $
 \widetilde{V}\geq 0. $

 Since $\underline{V}\,$ satisfies
 \bee
 \left\{
 \begin{array}{l}
 d\underline{V}\,_t=-\,[\,L\,\underline{V}\,_t+ M^k\,\underline{Z}\,^k_t
 +(\,g_t-L\,\underline{V}\,_t-M^k\,\underline{Z}\,_t^k\,\,)\,]\,dt
 +\underline{Z}\,^k_t\,dB^k_t,
 \vspace{2mm} \\
 \underline{V}\,^+_T\leq\varphi\leq\varphi^+\,,
 \end{array}
 \right.
 \eee
 then the comparison theorem for linear BSPDE in \cite{Du} implies that
 $ \widetilde{V}\geq \underline{V}\,.$

 Define
 $$
 \overline{V}\triangleq\widetilde{V}+(1+|\,x\,|^{d^*+1})^{-1}.
 $$
 Then  $\overline{V}\,\in\mathbb{H}^{2,\,2},\,
 \overline{V}\,>\underline{V}\,^+,\;
 \overline{V}\,_T>\varphi\geq \underline{V}\,_T$,  and
 $$
 d\,\overline{V}\,_t=-\overline{g}\,_t\,dt+\overline{Z}\,^k_t\,dB^k_t$$
 with $
 \overline{g}\,=L\,\widetilde{V}_t
 + M^k\,\widetilde{Z}^k_t+\widetilde{f}_t
 \in\mathbb{H}^{0,\,2}$ and $\overline{Z}=\widetilde{Z}\in\mathbb{H}^{1,\,2}.
 $
 Hence, $\underline{V}\,,\,\overline{V},$ and $\varphi$ satisfy Assumptions V3$'$
 and V4. The estimate (\ref{eq7.1}) follows from Lemma~\ref{lem2.2}.

 In the following, we prove that Problems $\mathscr{O}_{tx}$ and $\sD_{tx}$ are equivalent.
 We firstly claim
\be                                                                                          \label{eq7.3}
 \mathbb{E}\,\Big[\,R_t(x;\tau_1,\,\tau_2)\,\Big|\,{\cal F}_t\,\Big]
 \geq\mathbb{E}\,\Big[\,P_t(x;\tau_1)\,\Big|\,{\cal F}_t\,\Big],\quad
 \forall\;\;\tau_1,\,\tau_2\in{\cal U}\,_{t,\,T}.
\ee

 In fact, on the event of $\{\tau_1<\tau_2\}$, it is clear that
 $$
 P_t(x;\tau_1)=R_t(x;\tau_1,\tau_2).
 $$
 On the event of $\{\tau_1\geq\tau_2\}$, applying Theorem~\ref{thm3.1} and
 repeating the method in the proof of Theorem~\ref{thm4.2},  we deduce that
 \bee
 R_t(x;\tau_1,\tau_2)
 &\!\!\!=\!\!\!&\int_t^{\tau_2}f_u(X^{t,x}_u)\,du
 +\varphi(X^{t,x}_T)\,\chi_{\{\tau_2\geq T\}}
 +\overline{V}\,_{\tau_2}(X^{t,x}_{\tau_2})\,\chi_{\{\tau_2<T\}}
 \\[2mm]
 &\!\!\!\geq\!\!\!&\int_t^{\tau_2}f_u(X^{t,x}_u)\,du
 +\varphi(X^{t,x}_T)\,\chi_{\{\tau_2\geq T\}}
 +\widetilde{V}\,_{\tau_2}(X^{t,x}_{\tau_2})\,\chi_{\{\tau_2<T\}}
 \\[2mm]
 &\!\!\!=\!\!\!&\int_t^{\tau_2}f_u(X^{t,x}_u)\,du
 +\varphi(X^{t,x}_T)\,\chi_{\{\tau_2\geq T\}}
 +\widetilde{V}\,_{\tau_1}(X^{t,x}_{\tau_1})\,\chi_{\{\tau_2<T\}}
 +\int^{\tau_1}_{\tau_2}\widetilde{f}_u(X^{t,x}_u)\,du\qquad
 \\[2mm]
 &\!\!\!\!\!\!&-\int^{\tau_1}_{\tau_2}
 (\widetilde{Z}\,^k_u+M^k\widetilde{V}\,_u)(X^{t,x}_u)\,dB^k_u
 -\int^{\tau_1}_{\tau_2} (N^l\widetilde{V}\,_u)\,(X^{t,x}_u)\,dW^l_u
 \\[2mm]
 &\!\!\!\geq\!\!\!&P_t(x;\tau_1)
 -\int^{\tau_1}_{\tau_2} (\widetilde{Z}\,^k_u
 +M^k\widetilde{V}\,_u)(X^{t,x}_u)\,dB^k_u
 -\int^{\tau_1}_{\tau_2} (N^l\widetilde{V}\,_u)\,(X^{t,x}_u)\,dW^l_u.
 \eee
 Hence, we obtain
 $$
 R_t(x;\tau_1,\tau_2)
 \geq P_t(x;\tau_1)
 -\int^{\tau_1\vee\tau_2}_{\tau_2} (\widetilde{Z}\,^k_u
 +M^k\widetilde{V}\,_u)(X^{t,x}_u)\,dB^k_u
 -\int^{\tau_1\vee\tau_2}_{\tau_2} (N^l\widetilde{V}\,_u)\,(X^{t,x}_u)\,dW^l_u.
 $$
 Taking the condition expectation in the above inequality, we
 have (\ref{eq7.3}).

 If Problem $\sD_{tx}$ has a saddle point $(\tau_1^*,\tau_2^*)$, then
 we have that
 \bee
 \mathbb{E}\,\Big[\,P_t(x;\tau_1^*)\,\Big|\,{\cal F}_t\,\Big]
 &=&\mathbb{E}\,\Big[\,R_t(x;\tau_1^*,\,T)\,\Big|\,{\cal F}_t\,\Big]
 \geq \mathbb{E}\,\Big[\,R_t(x;\tau_1^*,\,\tau_2^*)\,\Big|\,{\cal F}_t\,\Big]
 \\[2mm]
 &\geq& \mathbb{E}\,\Big[\,R_t(x;\tau_1,\,\tau_2^*)\,\Big|\,{\cal F}_t\,\Big]
 \geq\mathbb{E}\,\Big[\,P_t(x;\tau_1)\,\Big|\,{\cal F}_t\,\Big],
 \eee
 where $\tau_1$ is an arbitrary stopping time in ${\cal U}\,_{t,\,T}$£¬
 and we have used (\ref{eq7.3}) in the last inequality. Hence, Problem~$\mathscr{O}_{tx}$ has an optimal stopping time $\tau_1^*\in{\cal U}\,_{t,\,T}$.

 Suppose that Problem $\mathscr{O}_{tx}$ has an optimal stopping time
 $\tau_1^*\in{\cal U}\,_{t,\,T}$. Then we choose $\tau^*_2=T$.
 We see that for any $\tau_1\in{\cal U}\,_{t,\,T}$,
 $$
 R_t(x;\tau_1,\,\tau^*_2)=P_t(x;\tau_1)
 $$
 and
 $$
 \mathbb{E}\,\Big[\,R_t(x;\tau^*_1,\tau_2^*)\,\Big|\,{\cal F}_t\,\Big]
 =\mathbb{E}\,\Big[\,P_t(x;\tau^*_1)\,\Big|\,{\cal F}_t\,\Big]
 \geq\mathbb{E}\,\Big[\,P_t(x;\tau_1)\,\Big|\,{\cal F}_t\,\Big]
 =\mathbb{E}\,\Big[\,R_t(x;\tau_1,\tau_2^*)\,\Big|\,{\cal F}_t\,\Big].
 $$
 On the other hand, according to  (\ref{eq7.3}), we have that
 $$
 \mathbb{E}\,\Big[\,R_t(x;\tau_1^*,\,\tau_2)\,\Big|\,{\cal F}_t\,\Big]
 \geq\mathbb{E}\,\Big[\,P_t(x;\tau_1^*)\,\Big|\,{\cal F}_t\,\Big]
 =\mathbb{E}\,\Big[\,R_t(x;\tau^*_1,\tau_2^*)\,\Big|\,{\cal F}_t\,\Big].
 $$
 Hence, $(\tau_1^*,\tau_2^*)$
 is a saddle point of Problem $\sD_{tx}$. Until now, we have proved that Problems
$\mathscr{O}_{tx}$ and $\sD_{tx}$ are equivalent.

 Let Assumptions V1, V2, O1, and O2 are satisfied.
 Denote by $(\widehat{V},\,\widehat{Z})$ the solution of the
 following BSPDE:
 \be                                                                                          \label{eq7.4}
 \left\{
 \begin{array}{l}
 d\widehat{V}_t=-({\cal L}\,\widehat{V}_t
 + {\cal M}^k\,\widehat{Z}^k_t+\widehat{f}_t)\,dt
 + \widehat{Z}^k_t\,dB^k_t;
 \vspace{2mm} \\
 \widehat{V}_T(x)=\varphi^+\,,
 \end{array}
 \right.
\ee
 where ${\cal L}$ and ${\cal M}$ are defined in (\ref{eq1.3}) and $\widehat{f}$
 is defined as
 $$
 \widehat{f}=\max\{\,f,\,0,\,g-{\cal L}\,\underline{V}\,
 -{\cal M}^k\,\underline{Z}\,^k\,\}.
 $$
 Moreover, we define
 $$
 \overline{V}\,=\widehat{V}+(1+|\,x\,|^{d^*+1})^{-1}.
 $$

 Repeating the above argument,  we derive that BSPDE (\ref{eq7.4}) has a strong
 solution $(\widehat{V},\,\widehat{Z})\in \mathbb{H}^{2,\,2}\times
 \mathbb{H}^{1,\,2}$. Moreover, we have $\widehat{V}\geq \underline{V}\,^+,\,
 \overline{V}>\underline{V}\,$ and
 $\underline{V}\,,\,\overline{V},$ and $\varphi$ satisfy Assumptions V3$'$ and V4.
The estimate (\ref{eq7.1}) follows from Lemma~\ref{lem2.2}.\smallskip

 So, BSPDVI~\eqref{BSPDI2} has a unique strong solution $(V,Z,k^+,k^-)$ by
 Theorem~\ref{thm5.3}. \smallskip

 On the other hand, since $\underline{V}\,\leq \widehat{V}<\overline{V}\,$
 and $(\,\widehat{V},\widehat{Z}\,)$ is the strong solution of BSPDE (\ref{eq7.4}),
 then $(\widehat{V},\,\widehat{Z},\,0,0)$ is
 the strong solution of the following BSPDVI:
 \bee
 \left\{
 \begin{array}{l}
 d\widehat{V}_t=-({\cal L} \widehat{V}_t+{\cal M}^k \widehat{Z}^k_t
 +\widetilde{f}_t)\,dt +\widehat{Z}^k_t\,dB^k_t,\qquad
 \mbox{if}\;\;\underline{V}\,_t<\widehat{V}_t<\overline{V}\,_t\,;
 \vspace{2mm} \\
 d\widehat{V}_t\leq-({\cal L} \widehat{V}_t+{\cal M}^k \widehat{Z}^k_t
 +\widetilde{f}_t)\,dt+\widehat{Z}^k_t\,dB^k_t,\qquad
 \mbox{if}\;\;\widehat{V}_t=\underline{V}\,_t\,;
 \vspace{2mm} \\
 d\widehat{V}_t\geq-({\cal L} \widehat{V}_t+{\cal M}^k \widehat{Z}^k_t
 +\widetilde{f}_t)\,dt+\widehat{Z}^k_t\,dB^k_t,\qquad
 \mbox{if}\;\;\widehat{V}_t=\overline{V}\,_t\,;
 \vspace{2mm} \\
 \widehat{V}_T(x)=\varphi^+(x)\,.
 \end{array}
 \right.
 \eee
 In view of Theorem~\ref{thm5.2}, $\widehat{V}\geq V$ and
 $\overline{V}\,>V$. So, we deduce that $k^-=0$ a.e. in
 $\Omega\times Q$ and $(V,Z,k^+)$ is the strong solution of BSPDVI~\eqref{BSPDVI1}.\smallskip

 On the other hand, in view of Theorem~\ref{thm7.1},  the strong solution
 of BSPDVI~\eqref{BSPDVI1} is unique. So, the unique strong solutions of
 BSPDVI~\eqref{BSPDVI1} and ~\eqref{BSPDI2} coincide.
 \hfill$\Box$\medskip

 Recalling Lemma~\ref{lem7.2} and Remark~\ref{rem4.1} in Section 4,
 and Theorem~\ref{thm5.3} in Sections 5, we have

\begin{thm}                                                                                       \label{thm7.3}
 Let Assumptions V1, V2, O1, and O2 be satisfied.
 Then BSPDVI~(\ref{BSPDVI1}) has a unique strong solution $(V,Z,k^+)$ such that
 \bee
 &\!\!\!\!\!\!&\|V\|_{2,\,2}+ |\|V|\|_{1,\,2}+\|Z\|_{1,\,2}
 +\|k^+\|_{0,\,2}
 \\[2mm]
 &\!\!\!\leq\!\!\!& C(\kappa,K,T)\,\Big(\,
 \mathbb{E}\,[\,|\,\varphi\,|\,_{1,\,2}\,]
 +\|f\|_{0,\,2}+\|\,\underline{V}\,\|_{2,\,2}
 +\|\,\underline{Z}\,\|_{1,\,2}+\|\,\underline{g}\,\|_{0,\,2}\,\Big).\qquad
 \eee
 Moreover,  the strong solution of  BSPDVI~\eqref{BSPDVI1}
 coincides with the value of Problem $\mathscr{O}$ if (\ref{eq4.4}) holds.
\end{thm}

 Identically as in the proof of Theorem~\ref{thm6.1}, we have

\begin{thm}                                                                                      \label{thm7.4}
 Let Assumptions V1, V2, V5(i), O1, and O2 be satisfied, and  the functions $\underline{f}\,$ and
 $\,\underline{\varphi}\,$ be increasing (resp. decreasing) in
 $x_i$, with $\underline{f}$ and $\underline{\varphi}$ be defined in
 Assumption V6(i).  Then assertions (i) and (ii) in Theorem~\ref{thm6.1} hold.
 \smallskip

 Moreover, if Assumption V5(j) with $i\neq j$ is satisfied, and
 $\underline{f}$ and $\underline{\varphi}\,$ are monotone in
 $x_j$, then the free boundary $\overline{S}\,_i$
 is monotone in $x_j$  for any
 $x_{\!\not{\, i,}\not{\, j}}\in\mathbb{R}^{d^*-2}$, a.e. in $\Omega\times[\,0,T\,]$.
 \end{thm}

\end{document}